\documentclass[a4paper, 10pt, leqno]{amsbook} 

\usepackage[T1]{fontenc}
\usepackage[latin1]{inputenc}
\usepackage{ngerman}
\usepackage{amsmath}
\usepackage{amssymb}
\usepackage{setspace}

\usepackage{color}

\newenvironment{beweis}{\textsc{Proof:}\;\;}{\rightline{$\Box$}\\}

\begin{document}
\setcounter{page}{1}
\newtheorem{theorem}{\textbf{Theorem}}[chapter]
\newtheorem{proposition}[theorem]{\textbf{Proposition}}
\newtheorem{lemma}[theorem]{\textbf{Lemma}}
\newtheorem{remark}[theorem]{\textbf{Remark}}
\newtheorem{example}[theorem]{\textbf{Example}}
\newtheorem{corollary}[theorem]{\textbf{Corollary}}
\newtheorem{definition}[theorem]{\textbf{Definition}}

\vspace{1cm}

\centerline{\textsc{\Large\textbf{Gauge Deformations and}}}
\centerline{\textsc{\Large\textbf{Embedding Theorems for Special Geometries}}}

\vspace{3cm} \centerline{Sebastian Stock\footnote{sstock@math.uni-koeln.de}}\centerline{Department of Mathematics} \centerline{University of Cologne}

\vspace{2.5cm}

\centerline{\textsc{\textbf{Abstract}}}
\begin{center}
\parbox{12cm}{We reduce the embedding problem for hypo $SU(2)$ and $SU(3)$-structures to the embedding problem for hypo $G_{2}$-structures into parallel $\text{Spin(7)}$-manifolds. The latter will be described in terms of gauge deformations. This description involves the intrinsic torsion of the initial $G_{2}$-structure and allows us to prove that the evolution equations, for all of the above embedding problems, do not admit non-trivial longtime solutions.}
\end{center}

\vspace{1cm}
{\large \textbf{Acknowledgement}}

Part of this work was done during a stay at McMaster University in Hamilton, Ontario, Canada. I would like to thank especially Prof. McKenzie Wang for our meetings and many helpful discussions on the topic.

\vspace{2cm}
{\large \textbf{Introduction}}\\

In \cite{hit1} N.Hitchin introduced a flow equation for cocalibrated $G_{2}$-structures on a manifold $M$, whose solutions yield parallel $\text{Spin}(7)$-structures on $I\times M$, for some interval $I\subset\mathbb{R}$. In this sense, a solution of the flow equation embeds the initial $G_{2}$-structure into a manifold with a parallel $\text{Spin}(7)$-structure and is therefore called a solution of the embedding problem for the initial structure. Similar equations are known for embedding $SU(2)$-structures in dimension five and $SU(3)$-structures in dimension six into manifolds with a parallel $SU(3)$ and $G_{2}$-structure, respectively, cf. \cite{co},\cite{cosa},\cite{cor},\cite{fer1},\cite{fer2}. The natural candidates for solving the embedding problem are so-called hypo structures. In the Gray-Hervella classification these are the type of structures induced on hypersurfaces of spaces with a parallel structure. Hypo $SU(3)$-structures are also called half-flat structures, whereas hypo $G_{2}$-structures are often called cocalibrated structures. R.Bryant shows in \cite{bry} that in the real analytic category, the embedding problem for hypo $SU(3)$ and $G_{2}$-structures can be solved. Bryant also provided counterexamples in the smooth category. The embedding problem for $SU(2)$-structures in dimension five was solved by D. Conti and S. Salamon in \cite{cosa}, cf. also \cite{co}.

The purpose of this article is to describe a unifying approach to all of the above embedding problems. We reduce the $SU(2)$ and $SU(3)$ embedding problem to the $G_{2}$-case, which will be studied in terms of gauge deformations, i.e. automorphism of the tangent bundle. Since the structure tensor $\varphi\in\Omega^{3}(M)$ of a $G_{2}$-structure is stable, any smooth deformation $\varphi_{t}$ can be described by a family of gauge deformations $A_{t}\in C^{\infty}(\text{Aut}(TM))$ via $\varphi_{t}=A_{t}\varphi$. It seems to be coincidence, that in the $G_{2}$-case, the intrinsic torsion $\mathcal{T}$ takes values in the $G_{2}$-module $\mathfrak{gl}(7)$ and therefore can be regarded again as an (infinitesimal) gauge deformation.
In Proposition \ref{propeea} we show that the intrinsic torsion flow for $G_{2}$-structures
$$\dot{A}_{t}=\mathcal{T}_{t}\circ A_{t}$$
can be regarded as a generalization of Hitchin's flow equation, and hence as a generalization of the $SU(2)$, $SU(3)$ and $G_{2}$-embedding problem. We describe the evolution of the metric and the intrinsic torsion under the intrinsic torsion flow, cf. Theorem \ref{thmgt}. As a consequence of the Cheeger-Gromoll Splitting Theorem, we prove in Theorem \ref{thmltex} and Corollary \ref{corltex} that there are no nontrivial longtime solutions for the embedding problem.\\

In chapter 2 we develop a conservation law for certain integral curves in Fr\'{e}chet spaces, cf. Corollary \ref{corintcurv}. The basic idea stems from finite dimensional geometry: If a vector field $X$ is tangent to some submanifold $N$, then any integral curve of $X$, which lies initially in $N$, stays in $N$ for all times. This does not hold for arbitrary integral curves in Fr\'{e}chet spaces, but the Cauchy-Kowalevski Theorem states - beyond the existence - that the integral curves in question can be developed in a (convergent) power series. This property allows us to prove that the intrinsic torsion flow preserves certain compatibility conditions, which implies that for any real analytic hypo $SU(2)$, $SU(3)$ and $G_{2}$-structure on a compact manifold, the embedding problem admits a unique real analytic solution. Moreover, the solution can be described by a family of gauge deformations
$$A_{t}=\sum_{k=0}^{\infty}\frac{t^{k}}{k!}A^{(k)}_{0},$$
where the series converges in the $C^{\infty}$-topology on $C^{\infty}(\text{End}(TM))$.\\

Our technique seems to be applicable to a wide class of evolution problems, where the initial structure is real analytic. For instance, instead of embedding a certain $G$-structure into a manifold with a parallel structure, one can ask for an embedding into a space with a nearly parallel structure, cf. \cite{co}. The Cauchy-Kowalevski Theorem \ref{thmckglobal} ensures the existence of a solution for the corresponding evolution equations. This solution has to satisfy certain (non linear) compatibility conditions. Since Corollary \ref{corintcurv} can be generalized to integral curves in Fr\'{e}chet manifolds, it suffices to show that the evolution equations define a vector field which is tangent to the compatibility conditions.\\

\setcounter{chapter}{1}
\setcounter{theorem}{0}
\vspace{1cm}
\begin{center}{\large \textbf{1. The Embedding Problem for Special Geometries}}\end{center}
\vspace{1cm}

A $G$-structure on a manifold $M$ is a reduction of the structure group of the frame bundle to a certain Lie subgroup $G\subset GL(n)$. We are interested in the cases where $M=M^{n}$ is a compact oriented manifold of dimension $n\in\{5,6,7,8\}$ and
$$G\in\{SU(2),SU(3),G_{2},\text{Spin}(7)\}.$$
The above groups can be realized as the isotropy group of certain model forms $\varphi_{0}\in\Lambda^{k}\mathbb{R}^{n\ast}$, under the natural action of $GL^{+}(n)$. The corresponding forms $\varphi\in C^{\infty}(\Lambda^{k}T^{\ast}M)$ on $M$ are called the structure tensors of the $G$-structure. A positive basis of $T_{p}M$, for which $\varphi\cong\varphi_{0}$, is called a Cayley frame for $\varphi$ and we say that $\varphi$ is of type $\varphi_{0}$. Since $G\subset SO(n)$, the structure tensors induce a metric $g=g(\varphi)$ on $M$ and we denote by $\nabla^{g}$ the Levi-Civita connection of the metric. The structure is called parallel if $\nabla^{g}\varphi=0$ holds. In the above cases, $\nabla^{g}\varphi=0$ can be translated into the apparently weaker conditions $d\varphi=d\ast\varphi=0$.\\


\begin{example}\label{exsp7}\textup{
A $\text{Spin}(7)$-structure on $M^{8}$ can be described by a single $4$-form $\Psi$ of type
\begin{equation*}
\begin{split}
\Psi_{0}&=e^{3456}+e^{3478}+e^{5678}-e^{2358}+e^{2468}-e^{2457}-e^{2367}\\
&\;\;\;+e^{1357}-e^{1467}-e^{1458}-e^{1368}+e^{1234}+e^{1256}+e^{1278},\\
\end{split}
\end{equation*}
where $\text{Iso}_{GL(8)}(\Psi_{0})=\text{Spin}(7)$ holds. A $\text{Spin}(7)$-structure is parallel if $d\Psi=0$ holds, cf. \cite{sal}.\\}
\end{example}


\begin{example}\label{exg2}\textup{
A $G_{2}$-structure on $M^{7}$ can be described by a single $4$-form $\psi$ of type
$$\psi_{0}=e^{2345}+e^{2367}+e^{4567}-e^{1247}+e^{1357}-e^{1346}-e^{1256},$$
where $\text{Iso}_{GL^{+}(7)}(\psi_{0})=G_{2}$ holds. Given an orientation $[\varepsilon]$ for $M^{7}$, we can define a positive volume element $\varepsilon:=\varepsilon(\psi)\in\Lambda^{7}T^{\ast}M^{7}$ and a metric $g=g(\psi)$, cf. \cite{hit1}. Then the Hodge dual $\varphi:=\ast_{\psi}\psi$ is of model type
$$\varphi_{0}=e^{246}-e^{356}-e^{347}-e^{257}+e^{123}+e^{145}+e^{167}.$$
A $G_{2}$-structure is parallel if $d\varphi=d\psi=0$ holds.\\}
\end{example}


\begin{example}\label{exsu3}\textup{
A $SU(3)$-structure on $M^{6}$ can be described by a $4$-form $\sigma$ and a $3$-form $\rho$ of type
\begin{equation*}
\begin{split}
\sigma_{0}&=e^{1234}+e^{1256}+e^{3456},\\
\rho_{0}&=e^{135}-e^{245}-e^{236}-e^{146},\\
\end{split}
\end{equation*}
where $\text{Iso}_{GL^{+}(6)}(\sigma_{0},\rho_{0})=SU(3)$ holds. Given an orientation for $M^{6}$, we can define positive volume elements $\varepsilon:=\varepsilon(\sigma)=\varepsilon(\rho)\in\Lambda^{6}T^{\ast}M^{6}$, cf. \cite{hit1}. We consider $\sigma$ as an element $\sigma\in\text{Hom}(\Lambda^{2}TM^{6},\Lambda^{2}T^{\ast}M^{6})$ and define
$$\omega:=\frac{1}{2}\sigma(\omega^{\ast})\in\Lambda^{2}T^{\ast}M^{6},$$
where $\omega^{\ast}\in\Lambda^{2}TM^{6}$ is defined by $\sigma=\omega^{\ast}\otimes\varepsilon(\sigma)\in\Lambda^{2}TM^{6}\otimes\Lambda^{6}T^{\ast}M^{6}$. Then $\omega$ is of type $\omega_{0}=e^{12}+e^{34}+e^{56}$ and
\begin{equation*}
\begin{split}
2\alpha(I(\rho)X)\varepsilon&:=\rho\wedge(X\lrcorner\rho)\wedge\alpha,\\
\widehat{\rho}&:=-I(\rho)\lrcorner\rho,\\
2g(X,Y)\varepsilon&:=(X\lrcorner\rho)\wedge(Y\lrcorner\rho)\wedge\omega,\\
\end{split}
\end{equation*}
($X,Y\in TM^{6}$, $\alpha\in\Lambda^{1}T^{\ast}M^{6}$) define tensors of type $I_{0}=e^{1}\wedge e_{2}+..+e^{5}\wedge e_{6}$, $\widehat{\rho}_{0}=e^{136}-e^{246}+e^{235}+e^{145}$ and $g_{0}=\sum_{i=1}^{6}e^{i}\otimes e^{i}$, respectively. A $SU(3)$-structure is parallel if $d\omega=d\rho=d\widehat{\rho}=0$ holds.\\}
\end{example}


\begin{example}\label{exsu2}\textup{
A $SU(2)$-structure on $M^{5}$ can be described by a $2$-form $\omega_{1}$ and two $3$-forms $\rho_{2}$ and $\rho_{3}$ of type
\begin{equation*}
\begin{split}
\omega_{1}&=e^{23}+e^{45},\\
\rho_{2}&=e^{124}-e^{135},\\
\rho_{3}&=e^{125}+e^{134},\\
\end{split}
\end{equation*}
where $\text{Iso}_{GL^{+}(5)}(\omega_{1},\rho_{2},\rho_{3})=SU(2)$ holds, cf. Lemma \ref{lemdefsu2}. Given an orientation for $M^{5}$, we can define a positive volume element $\varepsilon:=\varepsilon(\omega_{1},\rho_{2},\rho_{3})\in\Lambda^{5}T^{\ast}M^{5}$, see Lemma \ref{lemepssu2}. Then
\begin{equation*}
\begin{split}
2\alpha(X)\varepsilon&:=(X\lrcorner\rho_{2})\wedge\rho_{2},\\
\omega_{2}(X,Y)\varepsilon&:=-(X\lrcorner\omega_{1})\wedge(Y\lrcorner\omega_{1})\wedge\rho_{2},\\
\omega_{3}(X,Y)\varepsilon&:=-(X\lrcorner\omega_{1})\wedge(Y\lrcorner\omega_{1})\wedge\rho_{3},\\
g(X,Y)\varepsilon&:=\alpha(X)\alpha(Y)\varepsilon+\alpha\wedge\omega_{1}\wedge(X\lrcorner\omega_{2})\wedge(Y\lrcorner\omega_{3}),\\
\end{split}
\end{equation*}
($X,Y\in TM^{6}$) define tensors of type $\alpha_{0}=e^{1}$, $\omega_{2}=e^{24}-e^{35}$, $\omega_{3}=e^{25}+e^{34}$ and $g_{0}=\sum_{i=1}^{5}e^{i}\otimes e^{i}$, respectively.\\}
\end{example}


In the previous examples, the model tensors in dimension $n+1$ can be constructed from the model tensors in dimension $n$. This is due to the fact that the inclusions
$$SU(2)\subset SU(3)\subset G_{2}\subset \text{Spin}(7)$$
can be realized as isotropy groups of certain unit vectors. In the following we will consider families of structures on $M$ which depend on a parameter $t\in I\subset\mathbb{R}$ and evolve under certain evolution equations. These equations actually guarantee that the induced structure on $I\times M$ is parallel. For instance, consider a family of $G_{2}$-structures $\psi_{t}$ on $M^{7}$, $t\in I$. Then
$$\Psi:=\psi_{t}+dt\wedge\varphi_{t}$$
defines a $\text{Spin}(7)$-structure on $M^{8}:=I\times M^{7}$ and
$$d^{8}\Psi=d^{7}\psi_{t}+dt\wedge\dot{\psi}_{t}-dt\wedge d^{7}\varphi_{t}=d^{7}\psi_{t}+dt\wedge(\dot{\psi}_{t}-d^{7}\varphi_{t}),$$
where $d^{7}$, $d^{8}$ denotes the exterior differential on $M^{7}$, $M^{8}$, respectively. Hence the $\text{Spin}(7)$-structure is parallel if and only if $d^{7}\psi_{t}=0$ and $\dot{\psi}_{t}=d\varphi_{t}$. The second equation can be regarded as an evolution equation for the initial structure $\varphi:=\varphi_{t=0}$, whereas $G_{2}$-structures with $d\psi_{t}=0$ are called hypo structures. Note that the evolution equation preserves the hypo condition $d\psi=0$. In the following Proposition we list the lifting maps for the $SU(2)$, $SU(3)$ and $G_{2}$-case, the hypo condition for the initial structure and the evolution equations to obtain parallel structures on $I\times M^{n}$.\\


\begin{proposition}\label{propee}\textup{
Let $M^{n}$ be a manifold of dimension $n\in\{5,6,7\}$, equipped with a family of
$$G_{n}:=\begin{cases}\;SU(2) &,n=5\\
\;SU(3) &,n=6\\
\;G_{2} &,n=7\\
\big(\text{Spin}(7) &,n=8\big)\\
\end{cases}$$
structures. Then the lift in the following table defines a $G_{n+1}$-structure on $M^{n+1}:=I\times M^{n}$:\\
\begin{center}
\begin{tabular}{|l||l|l|l|}\hline
$n$ &Lift &Hypo Condition &Evolution\\ \hline \hline
5 &$\omega:=\omega_{1}+dt\wedge\alpha$ &$0=d\omega_{1}$ &$\dot{\omega}_{1}=d\alpha$\\
&$\sigma:=\frac{1}{2}\omega_{1}^{2}+dt\wedge\alpha\wedge\omega_{1}$ &$0=d\rho_{2}$ &$\dot{\rho}_{2}=d\omega_{3}$\\
&$\rho:=-\rho_{3}+dt\wedge\omega_{2}$ &$0=d\rho_{3}$ &$\dot{\rho_{3}}=-d\omega_{2}$\\
&$\widehat{\rho}:=\rho_{2}+dt\wedge\omega_{3}$ & &\\ \hline
6 &$\varphi:=\rho+dt\wedge\omega$ &$0=d\rho$ &$\dot{\rho}=d\omega$\\
&$\psi:=\sigma-dt\wedge\widehat{\rho}$ &$0=d\sigma$ &$\dot{\sigma}=-d\widehat{\rho}$\\ \hline
7 &$\Psi:=\psi+dt\wedge\varphi$ &$0=d\psi$ &$\dot{\psi}=d\varphi$\\ \hline
\end{tabular}
\end{center}
\textbf{(1)} The structure on $M^{n+1}$ is parallel if and only if the initial structure is hypo and evolves according to the evolution equations from the table.\\\\
\textbf{(2)} The metric of the $G_{n+1}$-structure on $I\times M^{n}$ is given by $g=dt^{2}+g_{t}$, where $g_{t}$ is the family of metrics induced by the $G_{n}$-structures on $M^{n}$.}
\end{proposition}

\begin{beweis}
Choosing a Cayley frame $(E_{1}(t),..,E_{n}(t))$ for the family of $G_{n}$-structures, we obtain a Cayley frame for the lift by
$$(\frac{d}{dt},E_{1}(t),..,E_{n}(t)).$$
This proves that the lift actually defines a $G_{n+1}$-structure and that the metric is given by the formula in (2). The proof of (1) is similar to the $G_{2}$-case.\\
\end{beweis}


\begin{definition}\label{defemprob}\textup{
Let $M^{n}$ be a manifold of dimension $n\in\{5,6,7\}$, equipped with a hypo $G_{n}$-structure. A family of $G_{n}$-structures which solves the evolution equations from Proposition \ref{propee} and equals the initial structure at $t=0$ is called a solution of the embedding problem for the initial $G_{n}$-structure.}\\
\end{definition}


The lift from Proposition \ref{propee} does not preserve the hypo condition. This motivates\\

\begin{definition}\label{defhypolift}\textup{
Let $M^{n}$ be a manifold of dimension $n\in\{5,6\}$, equipped with a $G_{n}$-structure. We call
\begin{center}
\begin{tabular}{|l|l|}\hline
$n=5$ &$n=6$\\ \hline \hline
$\omega:=\omega_{3}+d\theta\wedge\alpha$ &$\varphi:=-\widehat{\rho}+d\theta\wedge\omega$ \\
$\sigma:=\frac{1}{2}\omega_{3}^{2}+d\theta\wedge\rho_{3}$ &$\psi:=\sigma-d\theta\wedge\rho$\\
$\rho:=\rho_{2}-d\theta\wedge\omega_{1}$ &\\
$\widehat{\rho}:=-\alpha\wedge\omega_{1}-d\theta\wedge\omega_{2}$ &\\ \hline
\end{tabular}
\end{center}
the hypo lift of the $G_{n}$-structure to $S^{1}\times M^{n}$. Conversely, given a $G_{n+1}$-structure on a manifold $M^{n+1}$, we obtain a $G_{n}$-structure on any oriented hypersurface $i:M^{n}\hookrightarrow M^{n+1}$ by
\begin{center}
\begin{tabular}{|l|l|}\hline
$n=5$ &$n=6$ \\ \hline \hline
$\omega_{1}:=-i^{\ast}(\frac{\partial}{\partial\theta}\lrcorner\rho)$ &$\rho:=-i^{\ast}(\frac{\partial}{\partial\theta}\lrcorner\psi)$\\
$\rho_{2}:=i^{\ast}\rho$ &$\sigma:=i^{\ast}\psi$\\
$\rho_{3}:=i^{\ast}(\frac{\partial}{\partial\theta}\lrcorner\sigma)$ &\\ \hline
\end{tabular}
\end{center}
where $\frac{\partial}{\partial\theta}$ is a global vector field along $i:M^{n}\hookrightarrow M^{n+1}$, which is orthonormal to $M^{n}$. We call the $G_{n}$-structure the structure induced by the $G_{n+1}$-structure and $\frac{\partial}{\partial\theta}$.\\}
\end{definition}


Note that we just applied the lifts from Proposition \ref{propee} to the structures
$$(\alpha,\omega_{3},-\omega_{1},-\omega_{2})=A(\alpha,\omega_{1},\omega_{2},\omega_{3}),$$
respectively,
$$(\omega,-\widehat{\rho},\rho)=I(\omega,\rho,\widehat{\rho}),$$
where $A\in GL^{+}(5)$ is defines by
$$A(e_{1},..,e_{5}):=(e_{1},e_{3},e_{4},e_{2},e_{5}).$$\\


\begin{lemma}\label{lemhypolift}\textup{
The hypo lift maps hypo structures to hypo structures.}
\end{lemma}

\begin{beweis}
In the $SU(2)$-case, we obtain $d\rho=0$ if $d\omega_{1}=d\rho_{2}=0$. The compatibility condition $\omega_{3}^{2}=\omega_{1}^{2}$ and $d\rho_{3}=0$ imply $d\sigma=0$. For a hypo $SU(3)$-structure we obtain immediately $d\psi=d\sigma+d\theta\wedge d\rho=0$.\\
\end{beweis}


We will now study the compatibility of the hypo lift with the evolution equations from Proposition \ref{propee}.\\

\begin{lemma}\label{lemhlee}\textup{
\textbf{(1)} Suppose $\psi$ is a family of $G_{2}$-structures on $M^{7}=S^{1}\times M^{6}$ which is the hypo lift of some family of $SU(3)$-structure $(\rho,\sigma)$ on $M^{6}$. Then
$$\dot{\psi}=d\varphi\;\;\;\Leftrightarrow\;\;\;\left\{\begin{array}{l}\dot{\rho}=d\omega\\ \dot{\sigma}=-d\widehat{\rho}\end{array}\right.$$
\textbf{(2)} Suppose $(\rho,\sigma)$ is a family of $SU(3)$-structures on $M^{6}=S^{1}\times M^{5}$ which is the hypo lift of some family of $SU(2)$-structure $(\omega_{1},\rho_{2},\rho_{3})$ on $M^{5}$. Then
$$\left.\begin{array}{l}\dot{\rho}=d\omega\\ \dot{\sigma}=-d\widehat{\rho}\end{array}\right\}\;\;\;\Leftrightarrow\;\;\;
\left\{\begin{array}{l}\dot{\omega}_{1}=d\alpha\\ \dot{\rho}_{2}=d\omega_{3}\\ \dot{\rho}_{3}=-d\omega_{2}\\ (\frac{1}{2}\omega_{3}^{2})^{\cdot}=d(\alpha\wedge\omega_{1}) \end{array}\right.$$}
\end{lemma}

\begin{beweis}
By assumption we have $\psi=\sigma-d\theta\wedge\rho$ and $\varphi=-\widehat{\rho}+d\theta\wedge\omega$. Hence
$$\dot{\psi}=\dot{\sigma}-d\theta\wedge\dot{\rho}\;\;\;\text{ and }\;\;\;d\varphi=-d\widehat{\rho}-d\theta\wedge d\omega$$
and part (1) follows. Similarly for part (2),
\begin{equation*}
\begin{array}{ll}
\omega=\omega_{3}+d\theta\wedge\alpha,\;\;\;\;\; &\sigma=\frac{1}{2}\omega_{3}^{2}+d\theta\wedge\rho_{3},\\
\rho=\rho_{2}-d\theta\wedge\omega_{1},\;\;\;\;\; &\widehat{\rho}=-\alpha\wedge\omega_{1}-d\theta\wedge\omega_{2}\\
\end{array}
\end{equation*}
gives
\begin{equation*}
\begin{split}
\dot{\rho}&=\dot{\rho}_{2}-d\theta\wedge\dot{\omega}_{1},\\
d\omega&=d\omega_{3}-d\theta\wedge d\alpha,\\
\end{split}
\end{equation*}
and
\begin{equation*}
\begin{split}
\dot{\sigma}&=(\frac{1}{2}\omega_{3}^{2})^{\cdot}+d\theta\wedge\dot{\rho}_{3},\\
-d\widehat{\rho}&=d(\alpha\wedge\omega_{1})-d\theta\wedge d\omega_{2}.\\
\end{split}
\end{equation*}\\
\end{beweis}


\begin{lemma}\label{lemishl}\textup{
Let $\psi$ be a $G_{2}$-structure on $M^{7}$ with metric $g$.\\
\textbf{(1)} If $M^{7}=S^{1}\times M^{6}$, then $\psi$ is the hypo lift of some $SU(3)$-structure on $M^{6}$ if and only if
$$L_{\frac{\partial}{\partial\theta}}\psi=0,\;\;\;\frac{\partial}{\partial\theta}\bot_{g} TM^{6}\;\;\;\text{ and }\;\;\;
g(\frac{\partial}{\partial\theta},\frac{\partial}{\partial\theta})=1.$$\\
\textbf{(2)} If $M^{7}=S^{1}_{2}\times S^{1}_{1}\times M^{5}$, then $\psi$ is the hypo lift of some $SU(2)$-structure on $M^{5}$ if and only if
$$L_{\frac{\partial}{\partial\theta_{i}}}\psi=0,\;\;\;
\frac{\partial}{\partial\theta_{i}}\bot_{g} TM^{5}\;\;\;\text{ and }\;\;\;
g(\frac{\partial}{\partial\theta_{i}},\frac{\partial}{\partial\theta_{j}})=\delta_{ij},$$\\
for $i,j=1,2$.}
\end{lemma}

\begin{beweis}
If $\psi$ is the hypo lift of some $SU(2)$ or $SU(3)$-structure, we get $L_{\frac{\partial}{\partial\theta_{i}}}\psi=0$ and the orthogonality condition on the $S^{1}$-directions. Conversely, we define forms $\sigma$ and $\rho$ on $M^{7}$ by
$$\psi=\underbrace{\frac{\partial}{\partial\theta}\lrcorner(d\theta\wedge\psi)}_{=:\sigma}
+d\theta\wedge\underbrace{(\frac{\partial}{\partial\theta}\lrcorner\psi)}_{=:-\rho}.$$
Since $\frac{\partial}{\partial\theta}$ is orthonormal to $M^{6}$ and $G_{2}$ acts transitively on $S^{6}$, we can find a Caley frame for which $\sigma$ and $\rho$ are of model type. Hence $(\sigma,\rho)$ defines a $SU(3)$-structure on each hypersurface $\{e^{i\theta}\}\times M^{6}$. Since
$$0=L_{\frac{\partial}{\partial\theta}}\sigma-d\theta\wedge L_{\frac{\partial}{\partial\theta}}\rho$$
implies $L_{\frac{\partial}{\partial\theta}}\sigma=L_{\frac{\partial}{\partial\theta}}\rho=0$, we see that $\sigma$ and $\rho$ are actually constant along the flow of $\frac{\partial}{\partial\theta}$. Part (2) of the Lemma follows similarly, using that $G_{2}$ acts transitively on pairs of orthonormal vectors.\\
\end{beweis}


\setcounter{chapter}{2}
\setcounter{theorem}{0}
\vspace{1cm}
\begin{center}\textsc{\large \textbf{2. Integral Curves in Fr\'{e}chet Spaces and the Cauchy-Kowalevski Theorem}}\end{center}
\vspace{1cm}

Hamilton \cite{ham} gives an introduction to Fr\'{e}chet manifolds which goes far beyond of what we require for our purposes. Although Proposition \ref{propintcurv} and Corollary \ref{corintcurv} can be generalized to Fr\'{e}chet manifolds, we focus on Fr\'{e}chet spaces to keep the technical effort at a minimum.  

A locally convex topological vector space $\mathcal{F}$ is a vector space with a collection of seminorms $\{\|.\|\}_{n\in N}$, i.e. functions $\{\|.\|\}_{n}:\mathcal{F}\rightarrow\mathbb{R}$ which satisfy
$$\|f\|\geq 0,\;\;\;\;\;\;\;\;\|f+g\|\leq\|f\|+\|g\|\;\;\;\;\text{ and }\;\;\;\;\|\lambda f\|=|\lambda|\|f\|,$$
for all $f,g\in\mathcal{F}$ and scalars $\lambda$. Such a family defines a unique topology which is metrizable if and only if $N$ is countable. In this case the topology is characterized by the property
$$\lim_{k\rightarrow\infty}f_{k}=f\in \mathcal{F}\;\;\;\;\Leftrightarrow\;\;\;\;\lim_{k\rightarrow\infty}\|f_{k}-f\|_{n}=0\text{ for all }n\in N.$$
The topology is Hausdorff if and only if $\|f\|_{n}=0$ for all $n\in N$, implies that $f=0$. The space is sequentially complete if every Cauchy sequence converges, where $f_{k}$ is a Cauchy sequence if it is a Cauchy sequence for every seminorm $\|.\|_{n}$.\\


\begin{definition}\label{deffrechet}\textup{
A Fr\'{e}chet space is a locally convex topological vector space, which is in addition metrizable, Hausdorff and complete.\\}
\end{definition}


\begin{example}\label{example}\textup{
Suppose $F\rightarrow M$ is a vector bundle over a compact manifold $M$. Then the vector space
$$\mathcal{F}:=C^{\infty}(F)$$
of smooth sections of $F$ is a Fr\'{e}chet space, where the collection of seminorms
$$\|f\|_{n}:=\sum_{j=0}^{n}\sup_{p\in M}|(\nabla^{(j)}f)(p)|$$
can be defined after choosing Riemannian metrics and connections on $TM$ and $F$, cf. \cite{ham} Example 1.1.5. The induced topology is the $C^{\infty}$ topology on $\mathcal{F}$.\\
Given an open subset $U\subset F$, we consider the subset of all sections in $\mathcal{F}$, whose image lies in $U$,
$$\mathcal{U}:=\{f\in\mathcal{F}\mid f(M)\subset U\}.$$
For $f\in\mathcal{U}$ we can find $\varepsilon>0$ such that
$$f\in B_{\varepsilon}^{0}(f):=\{\tilde{f}\in\mathcal{F}\mid \|\tilde{f}-f\|_{0}<\varepsilon\}\subset\mathcal{U}.$$
Since $B_{\varepsilon}^{0}(f)\subset\mathcal{F}$ is open, $\mathcal{U}$ is an open subset of the Fr\'{e}chet space $\mathcal{F}$.\\}
\end{example}


Smooth maps between Fr\'{e}chet spaces can be defined as follows: Let $U\subset \mathcal{F}$ be an open subset of a Fr\'{e}chet space $\mathcal{F}$ and $P:\mathcal{U}\rightarrow \mathcal{E}$ a continuous and nonlinear map into another Fr\'{e}chet space $\mathcal{E}$. We say that $P$ is $C^{1}$ on $\mathcal{U}$ if for every $f\in\mathcal{U}$ and every $v\in\mathcal{F}$ the limit
$$DP(f)v:=\lim_{t\rightarrow 0}\frac{1}{t}(P(f+tv)-P(f))$$
exists and the map $DP:\mathcal{U}\times\mathcal{F}\rightarrow\mathcal{E}$ is continuous. Consequently, we say that $P$ is $C^{k}$ on $\mathcal{U}$ if $P$ is $C^{k-1}$ and the limit
\begin{equation*}
\begin{split}
D^{(k)}P(f)\{v_{1},..,v_{k}\}:&=\\
\lim_{t\rightarrow 0}\frac{1}{t}&\bigg(D^{(k-1)}P(f+tv_{n})\{v_{1},..,v_{k-1}\}-D^{(k-1)}P(f)\{v_{1},..,v_{k-1}\}\bigg)
\end{split}
\end{equation*}
exists for all $f\in\mathcal{U}$ and $v_{1},..,v_{k}\in\mathcal{F}$, and the map $D^{(k)}P:\mathcal{U}\times\mathcal{F}\times..\times\mathcal{F}\rightarrow\mathcal{E}$ is continuous. We call $P$ a smooth map on $\mathcal{U}$ if $P$ is $C^{k}$ for all $k\in\mathbb{N}$. We summarize Corollary 3.3.5 and Theorem 3.6.2 from \cite{ham} in the following\\


\begin{theorem}\label{thmfrecn}\textup{
\textbf{(1)} If $P:\mathcal{U}\subset\mathcal{F}\rightarrow\mathcal{E}$ is $C^{1}$ and $c(t)\in\mathcal{U}\subset\mathcal{F}$ is a parametrized $C^{1}$ curve, then $P\circ c(t)$ is a parametrized $C^{1}$ curve and
$$\frac{\partial}{\partial t}(P\circ c(t))=DP(c(t))\dot{c}(t).$$
\textbf{(2)} If $P:\mathcal{U}\subset\mathcal{F}\rightarrow \mathcal{E}$ is $C^{k}$, then for every $f\in\mathcal{U}$
$$D^{(k)}P(f)\{v_{1},..,v_{k}\}$$
is completely symmetric and linear separately in $v_{1},..,v_{k}\in\mathcal{F}$.\\}
\end{theorem}


In the following we will consider curves $c(t)\in\mathcal{F}$ in a Fr\'{e}chet space $\mathcal{F}$, which are integral curves of a vector field that is tangent to some subspace $\mathcal{E}\subset\mathcal{F}$. In finite dimension we would expect that any such integral curve with $c(0)\in\mathcal{E}$ actually stays in the subspace for all times. This conclusion fails for Fr\'{e}chet spaces, as was pointed out to us by Christian B\"{a}r: Consider $\mathcal{F}:=C^{\infty}[1,2]$ and $\mathcal{E}:=\{0\}\subset\mathcal{M}$. Then
$$c_{t}(x):=\begin{cases}
(4\pi t)^{-\frac{1}{2}}\,\exp(-\frac{x^{2}}{4t}), &\text{ for $t>0$}\\
0, &\text{ for $t\leq0$}
\end{cases}$$
solves $\dot{c}_{t}=\Delta c_{t}=\partial^{2}c_{t}/\partial x^{2}$ and hence defines an integral curve of the vector field $X(c):=\Delta c$. Although $X$ is tangent to $\mathcal{E}$, i.e. $X(0)=0$, and $c_{0}=0\in\mathcal{E}$, the curve doesn't stay in $\mathcal{E}$, since $c_{t}\neq 0$, for $t>0$. Note also that $t\mapsto c_{t}(x)$ is not real analytic in $t=0$.\\


\begin{proposition}\label{propintcurv}\textup{
Suppose $\mathcal{E}\subset\mathcal{F}$ is a closed subspace of the Fr\'{e}chet space $\mathcal{F}$ and that $X:\mathcal{U}\subset\mathcal{F}\rightarrow\mathcal{F}$ is a smooth map defined on some open subset $\mathcal{U}\subset \mathcal{F}$.
Let $f\in\mathcal{F}$ and assume that
$$X_{\mid \mathcal{U}\cap\mathcal{E}_{f}}:\mathcal{U}\cap\mathcal{E}_{f}\rightarrow\mathcal{E},$$
where $\mathcal{E}_{f}:=\{f\}+\mathcal{E}$. If a smooth curve $c:(-\varepsilon,\varepsilon)\rightarrow\mathcal{F}$ satisfies
$$c(0)\in\mathcal{U}\cap\mathcal{E}_{f}\;\;\;\;\text{ and }\;\;\;\;X\circ c(t)=\dot{c}(t),$$
where $\dot{c}:(-\varepsilon,\varepsilon)\rightarrow\mathcal{F}$ is the derivative of $c(t)$ by $t$, then for all $k\geq 1$
$$c^{(k)}(0)\in\mathcal{E},$$
where $c^{(k)}:(-\varepsilon,\varepsilon)\rightarrow\mathcal{F}$ is the $k^{th}$ derivative of $c(t)$ by $t$.}
\end{proposition}

\begin{beweis}
First we prove by induction on $k$ that the $k^{th}$ differential $D^{(k)}X$ of $X:\mathcal{F}\rightarrow\mathcal{F}$ satisfies
\setcounter{equation}{0}
\begin{equation}
D^{(k)}X_{\mid\mathcal{U}\cap\mathcal{E}_{f}\times\mathcal{E}\times..\times\mathcal{E}}:\mathcal{U}\cap\mathcal{E}_{f}\times\mathcal{E}\times..\times\mathcal{E}\rightarrow\mathcal{E}.
\end{equation}
For $k=0$ this is just the assumption $X_{\mid \mathcal{U}\cap\mathcal{E}_{f}}:\mathcal{U}\cap\mathcal{E}_{f}\rightarrow\mathcal{E}$. For $v_{0}\in\mathcal{U}\cap\mathcal{E}_{f}$ and $v_{1},..,v_{k+1}\in \mathcal{E}$ we have by definition
\begin{equation*}
\begin{split}
&\;\;\;\;\;D^{(k+1)}X(v_{0})\{v_{1},..,v_{k+1}\}\\
&=\lim_{s\rightarrow 0}\frac{1}{s}\underbrace{(D^{(k)}X(\underbrace{v_{0}+sv_{k+1}}_{\in\mathcal{U}\cap\mathcal{E}_{f}\text{ for $s$ small}})\{v_{1},..,v_{k}\}-D^{(k)}X(v_{0})\{v_{1},..,v_{k}\})}_{\in\mathcal{E}\text{ by induction hypothesis}}
\end{split}
\end{equation*}
and since $\mathcal{E}$ is closed, we conclude that (1) holds for $k+1$. Next we show that
for $k\geq 0$ and any choice of smooth curves $t\mapsto v_{0}(t)\in\mathcal{U}$ and $t\mapsto v_{1}(t),..,v_{k}(t)\in\mathcal{F}$
\begin{equation}
\begin{split}
&\;\;\;\;\;\frac{\partial}{\partial t} D^{(k)}X(v_{0}(t))\{v_{1}(t),..,v_{k}(t)\}=D^{(k+1)}X(v_{0}(t))\{v_{1}(t),..,v_{k}(t),\dot{v}_{0}(t)\}\\
&\;\;\;\;\;\;\;\;\;\;\;\;\;\;\;\;\;\;\;\;\;\;\;\;\;\;\;\;\;\;\;\;\;\;\;\;\;\;\;\;\;\;\;\;\;\;\;\;\;\;\;\;\;\; +\sum_{j=1}^{k}D^{(k)}X(v_{0}(t))\{v_{1}(t),..,\dot{v}_{j}(t),..,v_{k}(t)\}
\end{split}
\end{equation}
holds. Applying Theorem \ref{thmfrecn} (1) to the map
$D^{(k)}X:\mathcal{U}\times\mathcal{F}\times..\times\mathcal{F}\rightarrow\mathcal{F}$, we get
\begin{equation*}
\begin{split}
&\;\;\;\;\;\frac{\partial}{\partial t} D^{(k)}X(v_{0}(t))\{v_{1}(t),..,v_{k}(t)\}\\
&=D(D^{(k)}X)(v_{0}(t),..,v_{k}(t))\{\dot{v}_{0}(t),..,\dot{v}_{k}(t)\}\\
&=\lim_{s\rightarrow 0}\frac{1}{s}\bigg(D^{(k)}X(v_{0}(t)+s\dot{v}_{0}(t))\{v_{1}(t)+s\dot{v}_{1}(t),..,v_{k}(t)+s\dot{v}_{k}(t)\}\\
&\;\;\;\;\;\;\;\;\;\;\;\;\;\;-D^{(k)}X(v_{0}(t))\{v_{1}(t),..,v_{k}(t)\} \bigg)
\end{split}
\end{equation*}
and (2) follows, since $D^{(k)}X$ is linear in the arguments in $\{...\}$, cf. Theorem \ref{thmfrecn} (2). We will now show by induction on $k$ that
$c^{(k)}(0)\in\mathcal{E}$ holds. For $k=1$ we have $\dot{c}(0)=X\circ c(0)\in\mathcal{E}$ by assumption. Since $\dot{c}(t)=X\circ c(t)=D^{(0)}X(c(t))$ and $c(t)\in\mathcal{U}$ for sufficiently small $t$, we can apply (2) to see that $c^{(k+1)}(t)$, again for sufficiently small $t$, can be expressed as a linear combination of
$$D^{(j)}X(c(t))\{v_{1}(t),..,v_{j}(t)\},$$
where $j\in\{1,..,k+1\}$ and $v_{1}(t),..,v_{j}(t)\in\{c^{(l)}(t)\mid 1\leq l\leq k\}$. Since $c(0)\in\mathcal{U}\cap\mathcal{E}_{f}$, we get from $c^{(1)}(0),..,c^{(k)}(0)\in\mathcal{E}$ and (1)
$$D^{(j)}X(c(0))\{v_{1}(0),..,v_{j}(0)\}\in\mathcal{E}$$
and hence $c^{(k+1)}(0)\in\mathcal{E}$.\\
\end{beweis}


The following corollary can be regarded as a conservation law for certain integral curves in Fr\'{e}chet spaces.\\

\begin{corollary}\label{corintcurv}\textup{
If the curve $c:(-\varepsilon,\varepsilon)\rightarrow\mathcal{F}$ from Proposition \ref{propintcurv} satisfies for all $t\in(-\varepsilon,\varepsilon)$
$$c(t)=\sum_{k=0}^{\infty}\frac{t^{k}}{k!}c^{(k)}(0)\in\mathcal{F},$$
where the series converges w.r.t. the Fr\'{e}chet topology in $\mathcal{F}$, then
$$c(t)-c(0)\in\mathcal{E},$$
for all $t\in(-\varepsilon,\varepsilon)$.}
\end{corollary}

\begin{beweis}
From Proposition \ref{propintcurv} we get $c^{(k)}(0)\in\mathcal{E}$ for all $k\geq 1$ and hence
$$c(t)-c(0)=\sum_{k=1}^{\infty}\frac{t^{k}}{k!}c^{(k)}(0)\in\mathcal{E},$$
since $\mathcal{E}\subset\mathcal{F}$ is closed and the series converges in $\mathcal{F}$.\\
\end{beweis}


A formal power series in $X=(X_{1},..,X_{n})$ with coefficients in $\mathbb{R}$ is an expression of the form
$$S(X)=\sum_{p\in\mathbb{N}^{n}}a_{p}X^{p},$$
where $a_{p}\in\mathbb{R}$ and $X^{p}:=X_{1}^{p_{1}}\cdot..\cdot X_{n}^{p_{n}}$, for $p=(p_{1},..,p_{n})\in\mathbb{N}^{n}$.
Given a formal power series $S(X)$, we define
\begin{equation*}
\begin{split}
\Gamma&:=\{r=(r_{1},..,r_{n})\mid r_{i}\geq 0\text{ and }\sum_{p\in\mathbb{N}^{n}}|a_{p}|\;r^{p}<\infty\}
\end{split}
\end{equation*}
and denote by $\Delta$ the interior of $\Gamma$, called the domain of convergence of the series. Hence the series
$$S(x)=\sum_{p\in\mathbb{N}^{n}}a_{p}x^{p}$$
is for every $x=(x_{1},..,x_{n})\in\mathbb{R}^{n}$ with $|x|=(|x_{1}|,..,|x_{n}|)\in\Gamma$ absolute convergent. We recall the following result:\\


\begin{proposition}\label{propconvps}\textup{
Suppose $S(X)$ is a formal power series with domain of convergence $\Delta$. For $\bar{x}=(\bar{x}_{1},..,\bar{x}_{n})\in\mathbb{R}^{n}$ with $|\bar{x}|\in\Delta$ and $r_{1},..,r_{n}$ with $0<r_{i}<|\bar{x}_{i}|$, define
$$K:=\{(x_{1},..,x_{n})\in\mathbb{R}^{n}\mid |x_{i}|\leq r_{i}\}.$$
\textbf{(1)} For any subset $P\subset\mathbb{N}^{n}$, the series
$$S_{P}(x):=\sum_{p\in P}a_{p}x^{p}$$
converges absolutely for all $x\in K$. In particular, the series $S(x):=\sum_{p\in\mathbb{N}^{n}}a_{p}x^{p}$ converges absolutely for $x\in K$.\\
\textbf{(2)} Suppose that $P_{N}\subset \mathbb{N}^{n}$ is a family of subsets, $N\in\mathbb{N}$, such that
$\lim_{N\rightarrow\infty}P_{N}=\mathbb{N}^{n}$. Then
$$S_{N}(x):=\sum_{p\in P_{N}}a_{p}x^{p}$$
converges uniformly on $K$ to the function $S:K\rightarrow\mathbb{R}$, $x\mapsto S(x)$.}\\
\end{proposition}

\begin{beweis}
Since $|\bar{x}|\in\Delta$ we can find $C>0$ such that
$$|a_{p}\bar{x}^{p}|\leq C,\;\;\;\text{ for all }p\in\mathbb{N}^{n}.$$
Hence for $x\in K$
\begin{equation*}
\begin{split}
|a_{p}x^{p}|&=|a_{p}\,\bar{x}_{1}^{p_{1}}\cdot..\cdot\bar{x}_{n}^{p_{n}}|\frac{|x_{1}^{p_{1}}\cdot..\cdot x_{n}^{p_{n}}|}{|\bar{x}_{1}^{p_{1}}\cdot..\cdot \bar{x}_{n}^{p_{n}}|}\\
&\leq C\left(\frac{r_{1}}{|\bar{x}_{1}|}\right)^{p_{1}}\cdot..\cdot\left(\frac{r_{n}}{|\bar{x}_{n}|}\right)^{p_{n}}.
\end{split}
\end{equation*}
Since $r_{i}/|\bar{x}_{i}|<1$, we can apply the method of majorants to see that $S_{P}(x)$ converges absolutely for $x\in K$. To prove uniform convergence consider
\begin{equation*}
\begin{split}
\sup_{x\in K}|S(x)-S_{N}(x)|&=\sup_{x\in K}|\sum_{p\in\mathbb{N}^{n}\setminus P_{N}}a_{p}x^{p}|\\
&\leq C\sum_{p\in\mathbb{N}^{n}\setminus P_{N}} \left(\frac{r_{1}}{|\bar{x}_{1}|}\right)^{p_{1}}\cdot..\cdot\left(\frac{r_{n}}{|\bar{x}_{n}|}\right)^{p_{n}}\\
\end{split}
\end{equation*}
Given $\varepsilon>0$, we can choose $M$ large, so that $\sum_{p_{i}=M+1}^{\infty}\left(\frac{r_{i}}{|\bar{x}_{i}|}\right)^{p_{i}}\leq\frac{\varepsilon}{nCC_{i}}$,
for $i=1,..,n$, where
$$C_{i}:=\sum_{\underset{\in\mathbb{N}^{n-1}}{(p_{1}..\hat{p}_{i}..p_{n})}}\left(\frac{r_{1}}{|\bar{x}_{1}|}\right)^{p_{1}}\cdot.. \cdot\widehat{\left(\frac{r_{i}}{|\bar{x}_{i}|}\right)}^{p_{1}} \cdot..\cdot\left(\frac{r_{n}}{|\bar{x}_{n}|}\right)^{p_{n}}<\infty\;\;\;\text{(geometric series)}.$$
The notation $\widehat{.}$ means that the corresponding factor is omitted. Since $\lim_{N\rightarrow\infty}P_{N}=\mathbb{N}^{n}$, we can find $N=N(M)$, such that $\{0,..,M\}^{n}\subset P_{N}$. Hence
\begin{equation*}
\begin{split}
\sup_{x\in K}|S(x)-S_{N}(x)|&\leq C\sum_{p\in\mathbb{N}^{n}\setminus \{0..M\}^{n}} \left(\frac{r_{1}}{|\bar{x}_{1}|}\right)^{p_{1}}\cdot..\cdot\left(\frac{r_{n}}{|\bar{x}_{n}|}\right)^{p_{n}}\\
&\leq C\sum_{i=1}^{n}\sum_{p_{i}=M+1}^{\infty}C_{i}\left(\frac{r_{i}}{|\bar{x}_{i}|}\right)^{p_{i}}\\
&\leq\varepsilon.
\end{split}
\end{equation*}\\
\end{beweis}


\begin{definition}\label{defrealana}\textup{
Let $U\subset \mathbb{R}^{n}$ open and $x_{0}\in U$.\\
\textbf{(1)} A function $f:U\rightarrow\mathbb{R}$ is called real analytic in $x_{0}\in U$ if there exists a formal power series $S$ with
$$f(x)=S(x-x_{0}),$$
for all $x$ in a neighborhood of $x_{0}$. \\
\textbf{(2)} A function $f:U\rightarrow\mathbb{R}$ is called real analytic in $U$ if $f$ is real analytic for every $x_{0}\in U$.\\
\textbf{(3)} A function $F=(f_{1},..,f_{m}):U\rightarrow\mathbb{R}^{m}$ is called real analytic in $U$ if each component $f_{i}:U\rightarrow\mathbb{R}$ is real analytic in $U$.}\\
\end{definition}

Note that the coefficients of $S$ can be computed in terms of partial derivatives, which shows that $S$ is uniquely determined by the condition $f(x)=S(x-x_{0})$. Moreover we have the following basic properties, cf. \cite{car} p.123:\\


\begin{lemma}\label{lemrealana}\textup{
\begin{enumerate}
\item[(1)] If $f:U\rightarrow\mathbb{R}$ is real analytic in $x_{0}\in U$, then it is differentiable in a neighborhood of $x_{0}$ and the derivatives are again real analytic functions in $x_{0}\in U$.\\
\item[(2)] If $f$ and $g$ are real analytic in $x_{0}$, then the product $fg$ is real analytic in $x_{0}$.\\
\item[(3)] If $f:U\rightarrow\mathbb{R}$ is real analytic, then $1/f$ is real analytic in all points $x\in U$, where $f(x)\neq 0$.\\
\item[(4)] Compositions of real analytic functions are again real analytic.\\
\end{enumerate}}
\end{lemma}


A manifold $M$ is called real analytic if it admits an atlas with real analytic transition functions. Similarly to the smooth category one can define real analytic vector bundles over $M$.\\


In the following we will develop a global version of the Cauchy-Kowalevski Theorem, cf. \cite{bryeds}, III. Theorem 2.1:\\

\begin{theorem}\label{thmck}\textup{
Let $t$ be a coordinate on $\mathbb{R}$, $x=(x_{i})$ be coordinates on
$\mathbb{R}^{n}$, $y=(y_{j})$ be coordinates on $\mathbb{R}^{s}$ and let $z=(z_{i}^{j})$ be coordinates on $\mathbb{R}^{ns}$. Let
$D\subset \mathbb{R}\times\mathbb{R}^{n}\times\mathbb{R}^{s}\times\mathbb{R}^{ns}$ open, and let $G:D\rightarrow\mathbb{R}^{s}$ be a
real-analytic mapping. Let $D_{0}\subset \mathbb{R}^{n}$ be open and $f:D_{0}\rightarrow\mathbb{R}^{s}$ be a real-analytic mapping with
Jacobian $Df(x)\in\mathbb{R}^{ns}$, i.e. $z_{i}^{j}(Df(x))=\partial f^{j}(x)/\partial x_{i}$, so that
$\{(t_{0},x,f(x),Df(x))\mid x\in D_{0}\}\subset D$ for some $t_{0}\in\mathbb{R}$.\\
Then there exists an open neighborhood $D_{1}\subset \mathbb{R}\times D_{0}$ of $\{t_{0}\}\times D_{0}$ and a real-analytic mapping
$F:D_{1}\rightarrow\mathbb{R}^{s}$ which satisfies
$$\begin{cases}
\frac{\partial F}{\partial t}(t,x)&=G(t,x,F(t,x),\frac{\partial F}{\partial x}(t,x))\\
F(t_{0},x)&=f(x)\;\;\;\text{ for all }x\in D_{0}.
\end{cases}$$
$F$ is unique in the sense that any other real-analytic solution of the above initial value problem agrees with $F$ in some neighborhood
of $\{t_{0}\}\times D_{0}$.\\}
\end{theorem}


\begin{remark}\label{remck}\textup{
Since the solution $F=(f_{i},..,f_{s}):D_{1}\rightarrow\mathbb{R}^{s}$ from Theorem \ref{thmck} is real analytic, we can develop each component in a convergent power series around $(t_{0},x_{0})=(0,0)\in D_{1}$, i.e.
$$f_{i}(t,x)=\sum_{k=0}^{\infty}\bigg(\sum_{p\in\mathbb{N}^{n}}a_{ikp}x^{p}\bigg)t^{k}
=\sum_{k=0}^{\infty}\bigg(\frac{1}{k!}f_{i}^{(k)}(0,x)\bigg)t^{k}.$$
Applying Proposition \ref{propconvps} (2) with $P_{N}:=\{0,..,N\}\times\mathbb{N}^{n}$ shows that
$$f_{i}^{N}(t,x)=\sum_{k=0}^{N}\bigg(\sum_{p\in\mathbb{N}^{n}}a_{ikp}x^{p}\bigg)t^{k}
=\sum_{k=0}^{N}\frac{t^{k}}{k!}f_{i}^{(k)}(0,x)$$
converges locally uniformly to the function $f_{i}(t,x)$, for $N\rightarrow\infty$. The partial derivatives of a formal power series $S(X)$ are defined by,
$$\frac{\partial S}{\partial X_{i}}:=\sum_{p\in\mathbb{N}^{n}}p_{i}a_{p}X_{1}^{p_{1}}\cdot..X_{i}^{p_{i}-1}..\cdot X_{n}^{p_{n}}.$$
The formal power series $\frac{\partial S}{\partial X_{i}}$ has the same domain of convergence $\Delta$ as the formal power series $S$. Moreover, the function $\frac{\partial S}{\partial X_{i}}:\Delta\rightarrow\mathbb{R}$ is the partial derivative of the function $S:\Delta\rightarrow\mathbb{R}$ w.r.t. $x_{i}$, cf. Satz 3.2 in \cite{car}. Hence we can apply again Proposition \ref{propconvps} (2) to see that all partial derivatives of the function
$f_{i}^{N}(t,x)$ converge locally uniformly to the corresponding partial derivative of $f_{i}(t,x)$. In summary, the functions
$$F_{N}(t,x):=\sum_{k=0}^{N}\frac{t^{k}}{k!}F^{(k)}(0,x)$$
converge, as $N\rightarrow\infty$, locally in $C^{\infty}$-topology to the solution $F(t,x)$ from Theorem \ref{thmck}.\\}
\end{remark}


\begin{definition}\label{defradiffop}\textup{
Suppose $M$ is a real analytic manifold and $\pi:V\rightarrow M$ is a rank $s$ real analytic vector bundle. We call a map
$$X:C^{\infty}(V)\rightarrow C^{\infty}(V)$$
a real analytic first order differential operator if every point of $M$ has a neighborhood $U\subset M$, which is the domain of a real analytic chart $u:U\rightarrow\mathbb{R}^{n}$, and there exists a real analytic trivialization $(\pi,v):V_{\mid U}\cong U\times\mathbb{R}^{s}$, together with a real analytic function
$$G:D\subset\mathbb{R}^{n}\times\mathbb{R}^{s}\times\mathbb{R}^{ns}\rightarrow\mathbb{R}^{s},$$
such that for every local section $c:U\subset M\rightarrow V$
$$v(X\circ c)=G(u,v\circ c,\frac{\partial c_{i}}{\partial u_{j}})$$
holds, where $c_{i}$ is the $i^{th}$ component of $v\circ c:U\rightarrow\mathbb{R}^{s}$.}
\end{definition}


We can now prove the following global version of the Cauchy-Kowalevski Theorem,\\

\begin{theorem}\label{thmckglobal}\textup{
Suppose $\pi:V\rightarrow M$ is a real analytic rank $s$ vector bundle over a compact real analytic manifold $M$. Let $X:C^{\infty}(V)\rightarrow C^{\infty}(V)$ be a real analytic first order differential operator and let $c_{0}\in C^{\infty}(V)$ be a real analytic section. Then the initial value problem
$$\begin{cases}
\dot{c}(t)=X\circ c(t)\\
c(0)=c_{0}\\
\end{cases}$$
has a unique real analytic solution $c:(-\varepsilon,\varepsilon)\rightarrow C^{\infty}(V)$, i.e. $c:(-\varepsilon,\varepsilon)\times M\rightarrow V$ is real analytic. Moreover, the solution $c(t)$ satisfies
$$c(t)=\sum_{k=0}^{\infty}\frac{t^{k}}{k!}c^{(k)}_{0},$$
where the series converges in the $C^{\infty}$ topology on $C^{\infty}(V)$.}
\end{theorem}

\begin{beweis}
We will first show that we can find local sections $c_{t}:U\subset M\rightarrow V$, which solve the initial value problem locally. Secondly, we prove that the compactness of $M$ ensures the existence of a global solution. Eventually we will use the uniqueness part of the Cauchy-Kowalevski Theorem to prove the uniqueness statement of the Theorem.

By Definition \ref{defradiffop} we can find a real analytic chart $u:U\subset M\rightarrow\mathbb{R}^{n}$ and a trivialization $(\pi,v):V_{\mid U}\cong U\times\mathbb{R}^{s}$, such that for each local section $c:U\subset M\rightarrow V$
\setcounter{equation}{0}
\begin{equation}\label{ckg1}
v(X\circ c)=G(u,v\circ c,\frac{\partial c_{i}}{\partial u_{j}})
\end{equation}
holds, where $G:D\subset\mathbb{R}^{n}\times\mathbb{R}^{s}\times\mathbb{R}^{ns}\rightarrow\mathbb{R}^{s}$ is real analytic. The map
$$f:D_{0}:=u(U)\subset\mathbb{R}^{n}\rightarrow\mathbb{R}^{s}\;\;\;\text{ with }\;\;\;f(x):=v\circ c_{0}\circ u^{-1}(x)$$
is real analytic and hence we can find by the Cauchy-Kowalevski Theorem a real analytic solution $F:(-\varepsilon,\varepsilon)\times\widetilde{D}_{0}\rightarrow\mathbb{R}^{s}$ of
$$\begin{cases}
\frac{\partial F}{\partial t}(t,x)&=G(x,F(t,x),\frac{\partial F}{\partial x}(t,x))\\
F(t_{0},x)&=f(x)\;\;\;\text{ for all }x\in D_{0},
\end{cases}$$
where $\widetilde{D}_{0}\subset D_{0}$ is open. Let $\widetilde{U}:=u^{-1}(\widetilde{D}_{0})\subset U$ and define for $t\in(-\varepsilon,\varepsilon)$
\begin{equation}\label{ckg2}
c(t):\widetilde{U}\subset M\rightarrow V\;\;\;\text{ by }\;\;\;c(t,p):=v_{p}^{-1}\circ F(t,u(p)),
\end{equation}
where $v_{p}:V_{p}\cong\mathbb{R}^{s}$ is the isomorphism induced by the local trivialization $(\pi,v)$. By definition, the map $c:(-\varepsilon,\varepsilon)\times\widetilde{U}\subset M\rightarrow V$ is real analytic and satisfies
\begin{equation}\label{ckg3}
c(0,p)=v_{p}^{-1}\circ F(0,u(p))=v_{p}^{-1}\circ f(u(p))=c_{0}(p).
\end{equation}
Now we have for $i=1,..,s$ and $j=1,..,n$
\begin{equation}\label{ckgcom}
\begin{split}
\frac{\partial(v_{i}\circ c_{t})}{\partial u_{j}}(p)&=\left.\frac{\partial}{\partial u_{j}}\right|_{p}\cdot (v_{i}\circ c_{t})
=(u^{-1}_{\ast}\left.\frac{\partial }{\partial x_{j}}\right|_{u(p)})\cdot(v_{i}\circ c_{t})\\
&=\left.\frac{\partial}{\partial x_{j}}\right|_{u(p)}\cdot(v_{i}\circ c_{t}\circ u^{-1})
=\left.\frac{\partial}{\partial x_{j}}\right|_{u(p)}\cdot F_{i}(t,.)\\
&=\frac{\partial F_{i}}{\partial x_{j}}(t,u(p)).\\
\end{split}
\end{equation}
Since by definition $v\circ c(t,p)=F(t,u(p))$ holds, we get from (\ref{ckg1}), applied to $c_{t}$
\begin{equation*}
\begin{split}
\dot{c}(t,p)&=v_{p}^{-1}\circ G(u(p),F(t,u(p)),\frac{\partial F}{\partial x}(t,u(p)))\\
&=v_{p}^{-1}\circ G(u(p),v\circ c_{t}(p),\frac{\partial (v_{i}\circ c_{t})}{\partial u_{j}}(p))\\
&=v_{p}^{-1}\circ v_{p}(X\circ c(t,p)),\\
&=X\circ c(t,p).
\end{split}
\end{equation*}
i.e. $c_{t}$ is the desired local solution of the initial value problem. Moreover, we get by Remark \ref{remck}
\begin{equation*}
\begin{split}
c(t,p)&=v_{p}^{-1}\circ F(t,u(p))=v_{p}^{-1}\big( \lim_{N\rightarrow\infty}\sum_{k=0}^{N}\frac{t^{k}}{k!}F^{(k)}(0,u(p))\big)\\
&=\lim_{N\rightarrow\infty}\sum_{k=0}^{N}\frac{t^{k}}{k!}v_{p}^{-1}\circ F^{(k)}(0,u(p))
=\lim_{N\rightarrow\infty}\sum_{k=0}^{N}\frac{t^{k}}{k!}c^{(k)}(0,p),\\
\end{split}
\end{equation*}
i.e.
\begin{equation}\label{ckg4}
c_{t}=\sum_{k=0}^{\infty}\frac{t^{k}}{k!}c^{(k)}_{0},
\end{equation}
where the series converges locally in $C^{\infty}$ topology.

Suppose now we apply the above construction to obtain two local sections
$$c_{1}(t):U_{1}\subset M\rightarrow V\;\;\;\text{ and }\;\;\;c_{2}(t):U_{2}\subset M\rightarrow V,$$
where $t\in(-\varepsilon,\varepsilon)$, $\varepsilon:=\text{min}\{\varepsilon_{1},\varepsilon_{2}\}$ and $U_{1}\cap U_{2}\neq\emptyset$. Since $c_{1}$ and $c_{2}$ both solve the initial value problem
$$\begin{cases}
\dot{c}_{i}(t)=X\circ c_{i}(t)\\
c_{i}(0)=c_{0},\\
\end{cases}$$
$i=1,2$, we see that $c_{1}(0)=c_{2}(0)$ and $\dot{c}_{1}(0)=\dot{c}_{2}(0)$ on $U_{1}\cap U_{2}$. Differentiating the equation $\dot{c}_{1}(t)=X\circ c_{1}(t)$, shows that $c^{(k+1)}_{1}(t)$ can be expressed as a linear combination of
$$D^{(j)}X(c_{1}(t))\{v_{1}(t),..,v_{j}(t)\},$$
where $j\in\{1,..,k+1\}$ and $v_{1}(t),..,v_{j}(t)\in\{c^{(l)}_{1}(t)\mid 1\leq l\leq k\}$, cf. the proof of Proposition \ref{propintcurv}. Now we obtain by induction $c^{(k)}_{1}(0)=c^{(k)}_{2}(0)$ on $U_{1}\cap U_{2}$, for all $k\in\mathbb{N}$. Hence (\ref{ckg4}) implies $c_{1}(t)=c_{2}(t)$ on $U_{1}\cap U_{2}$. If $M$ is compact, we can cover $M$ by finitely many domains $U_{1},..,U_{N}$ of local sections $c_{i}(t):U_{i}\subset M\rightarrow V$, which yield a global section $c(t):M\rightarrow V$, where $t\in(-\varepsilon,\varepsilon)$ and $\varepsilon:=\text{min}\{\varepsilon_{1},..,\varepsilon_{N}\}$. From (4) we get
$$c(t)=\sum_{k=0}^{\infty}\frac{t^{k}}{k!}c^{(k)}_{0},$$
and since $M$ is compact, the series converges in $C^{\infty}$ topology.

To prove uniqueness, suppose that we have two real analytic solutions $c_{1},c_{2}:(-\varepsilon,\varepsilon)\times M\rightarrow V$ of the initial value problem. By (\ref{ckg1}) we have for $k=1,2$ and $x\in u(U)\subset\mathbb{R}^{n}$
$$v(X\circ c_{k}(t)\circ u^{-1}(x))=G(x,v\circ c_{k}(t)\circ u^{-1}(x),\frac{\partial c_{ki}(t)}{\partial u_{j}}\circ u^{-1}(x)).$$
Now $F_{k}(t,x):=v\circ c_{k}(t)\circ u^{-1}(x)$ satisfies
$$\frac{\partial F_{k}}{\partial t}(t,x)=v\circ \dot{c}_{k}(t)\circ u^{-1}(x)=v\circ X\circ c_{k}(t)\circ u^{-1}(x)$$
and by (\ref{ckgcom})
$$\frac{\partial c_{ki}(t)}{\partial u_{j}}\circ u^{-1}(x)=\frac{\partial(v_{i}\circ c_{k}(t))}{\partial u_{j}}(u^{-1}(x))=\frac{\partial F_{ki}}{\partial x_{j}}(t,x),$$\\
for $i=1,..,s$ and $j=1,..,n$. Hence we showed
$$\frac{\partial F_{k}}{\partial t}(t,x)=G(x,F_{k}(t,x),\frac{\partial F_{ki}}{\partial x_{j}}(t,x)).$$
Since $F_{1}$ and $F_{2}$ are both real analytic and satisfy
$$F_{1}(0,x)=v\circ c_{1}(0)\circ u^{-1}(x)=v\circ c_{0}\circ u^{-1}(x)=v\circ c_{2}(0)\circ u^{-1}(x)=F_{2}(0,x),$$
the uniqueness part of the Cauchy-Kowalevski Theorem yields $F_{1}(t,x)=F_{2}(t,x)$, i.e. $c_{1}(t)=c_{2}(t)$.\\
\end{beweis}

\setcounter{chapter}{3}
\setcounter{theorem}{0}
\vspace{1cm}
\begin{center}\label{mc}{\textsc{\textbf{3. The Model Case $G_{2}\subset\text{Spin}(7)$}}}\end{center}
\vspace{1cm}

Lemma \ref{lemhlee} and \ref{lemishl} motivate the conjecture that the embedding problem for hypo $SU(2)$ and $SU(3)$-structures might be reduced to the embedding problem for $G_{2}$-structures. The reduction to the $G_{2}$-case has the advantage that no compatibility conditions are involved. To solve the embedding problem for hypo structures we consequently focus on studying the evolution equation
$$\dot{\psi}_{t}=d\varphi_{t}$$
on a compact seven dimensional manifold $M$. We will describe the solution $\psi_{t}$ by a family of gauge deformations, i.e.
$$\psi_{t}=A_{t}\psi,$$
where $A_{t}\in C^{\infty}(\text{Aut}(TM))$. Since the orbit of the model tensor $\psi\in\Lambda^{4}\mathbb{R}^{7\ast}$ is open, it follows that any smooth deformation $\psi_{t}$ of the initial structure $\psi$ can be described in such a way. The evolution equation $\dot{\psi}_{t}=d\varphi_{t}$ can be translated into an equation for the family of gauge deformations. This description involves the intrinsic torsion of the $G_{2}$-structure $\psi_{t}$. The intrinsic torsion $\mathcal{T}\in\text{End}(TM)$ of a $G_{2}$-structure $\varphi$ is defined by
$$\nabla^{g}_{X}\varphi=-\mathcal{T}X\lrcorner\psi,$$
where we used that
$$\nabla^{g}\varphi\in T^{\ast}M\otimes\Lambda^{3}_{7}T^{\ast}M$$
and
$$\Lambda^{3}_{7}T^{\ast}M:=\{\alpha\in\Lambda^{3}T^{\ast}M\mid \alpha=X\lrcorner\psi,\,X\in TM\},$$
cf. for instance \cite{bry}. From $\psi=\ast\varphi$ it follows that $d\psi=2\text{pr}_{\Lambda^{2}}(\mathcal{T})\wedge\varphi$ holds. So hypo $G_{2}$-structures are characterized by $\mathcal{T}\in S^{2}(TM)$ w.r.t. the metric $g$.\\


\begin{proposition}\label{propeea}\textup{
Suppose $\psi_{t}=A_{t}\psi$ is a family of $G_{2}$-structures on $M^{7}$, described by a family of gauge deformations $A_{t}\in C^{\infty}(\text{Aut}(TM^{7}))$. If $\mathcal{T}_{t}$ is the intrinsic torsion of $\psi_{t}$, then
$$\dot{\psi}_{t}=d\varphi_{t}\;\;\;\Leftrightarrow\;\;\;D_{\psi_{t}}(\dot{A}_{t}\circ A^{-1}_{t})=D_{\psi_{t}}(\mathcal{T}_{t}),$$
where
$$D_{\psi_{t}}:\text{End}(TM)\rightarrow\Lambda^{4}T^{\ast}M\;\;\;\text{ is defined by }\;\;\;A\mapsto\left.\frac{d}{ds}\right|_{s=0}\exp(sA)\psi_{t}.$$}
\end{proposition}

\begin{beweis}
Since clearly $\dot{\psi}_{t}=D_{\psi_{t}}(\dot{A}_{t}A_{t}^{-1})$, it suffices to observe that
\begin{equation*}
\begin{split}
D_{\psi_{t}}(\mathcal{T}_{t})(X_{1},..,X_{4})&=-\sum_{i=1}^{4}\psi_{t}(X_{1},..,\mathcal{T}_{t}X_{j},..,X_{4})\\
&=\sum_{i=1}^{4}(-1)^{i}\psi_{t}(\mathcal{T}_{t}X_{j},X_{1},..,\hat{X}_{j},..,X_{4})\\
&=\sum_{i=1}^{4}(-1)^{i+1}(\nabla^{g_{t}}_{X_{j}}\varphi_{t})(X_{1},..,\hat{X}_{j},..,X_{4})\\
&=d\varphi_{t}(X_{1},..,X_{4})
\end{split}
\end{equation*}
holds.\\
\end{beweis}


We can now compute the evolution of the metric and the torsion endomorphism.\\

\begin{theorem}\label{thmgt}\textup{
Let $\psi_{t}$ be a family of hypo $G_{2}$-structures on $M^{7}$, which evolves under the flow $\dot{\psi}_{t}=d\varphi_{t}$. Then the evolution of the underlying metric $g_{t}$ and the torsion endomorphism $\mathcal{T}_{t}$ are given by
\begin{equation*}
\begin{split}
&\dot{g}_{t}(X,Y)=2g_{t}(\mathcal{T}_{t}X,Y),\\
&\dot{\mathcal{T}}_{t}X=\text{Ric}_{t}X-\text{tr}(\mathcal{T}_{t})\mathcal{T}_{t}X,
\end{split}
\end{equation*}
where $\text{Ric}_{t}=\text{Ric}(g_{t})$ is the Ricci tensor of the metric $g_{t}$.}
\end{theorem}

\begin{beweis}
Writing $\psi_{t}=A_{t}\psi$, Proposition \ref{propeea} yields
$D_{\psi_{t}}(\dot{A}_{t}\circ A^{-1}_{t})=D_{\psi_{t}}(\mathcal{T}_{t})$. Since the evolution $\dot{\psi}_{t}=d\varphi_{t}$ preserves the hypo condition $d\psi_{t}=0$, or equivalently $\mathcal{T}_{t}\in S^{2}$ w.r.t. $g_{t}$, we get
$$\text{pr}_{S^{2}}(\dot{A}_{t}\circ A^{-1}_{t})=\mathcal{T}_{t},$$
since $\text{ker}(D_{\psi_{t}})=\mathfrak{g}_{2}$. Then we compute for $g_{t}=A_{t}g$
$$\dot{g}_{t}(X,Y)=2g_{t}(\text{pr}_{S^{2}}(\dot{A}_{t}\circ A^{-1}_{t})X,Y)=2g_{t}(\mathcal{T}_{t}X,Y).$$
The metric $g=dt^{2}+g_{t}$ on $I\times M^{7}$ has holonomy contained in $\text{Spin}(7)$ and hence is Ricci flat. The Gauss equations and the Codazzi-Mainardi equations yield
\begin{equation*}
\begin{split}
\dot{g}_{t}(X,Y)&=2g_{t}(\mathcal{W}_{t}X,Y),\\
g_{t}(\dot{\mathcal{W}}_{t}X,Y)&=\text{ric}_{t}(X,Y)-\text{tr}(\mathcal{W}_{t})g_{t}(\mathcal{W}_{t}X,Y),\\
\end{split}
\end{equation*}
where $\mathcal{W}_{t}X:=\nabla^{g}_{\Phi_{t\ast}X}\frac{d}{dt}$ is the Weingarten map and $\Phi_{t}$ is the flow of the vector field $\frac{d}{dt}$, cf. for instance \cite{baer}. So $\mathcal{W}_{t}=\mathcal{T}_{t}$ and the Theorem follows.\\
\end{beweis}


We will now apply the Cheeger-Gromoll Splitting Theorem to prove that the flow $\dot{\psi}=d\varphi$ does not admit nontrivial longtime solutions.\\

\begin{theorem}\label{thmltex}\textup{
Suppose $\psi$ is a hypo $G_{2}$-structures on a compact manifold $M^{7}$. Then the flow $\dot{\psi}_{t}=d\varphi_{t}$ is defined for all times
$t\in\mathbb{R}$ if and only if the initial structure is already parallel.\\}
\end{theorem}

\begin{beweis}
The metric on the product $M^{8}:=\mathbb{R}\times M^{7}$ has holonomy contained in $\text{Spin}(7)$ and hence is Ricci flat. Since $g=dt^{2}+g_{t}$, the first factor actually defines a line. Now we can apply the Cheeger-Gromoll Splitting Theorem and see that $M^{8}$ splits as a Riemannian product. Note that the line, i.e. the first factor of $M^{8}$, is actually the one dimensional factor that splits off in the decomposition as a Riemannian product, cf. Lemma 6.86 in \cite{bes}. Hence $g_{t}=g_{0}$ is constant and Theorem \ref{thmgt} yields $\mathcal{T}_{t}=0$.\\
\end{beweis}


In Lemma \ref{lemhlee} (1) we showed that a longtime solution of the $SU(3)$ embedding problem would yield a longtime solution for the $G_{2}$ embedding problem. Combining part (1) and (2) of Lemma \ref{lemhlee}, shows that a longtime solution of the $SU(2)$ embedding problem would also yield  a longtime solution for the $G_{2}$ embedding problem if in addition the equation $(\frac{1}{2}\omega_{3}^{2})^{\cdot}=d(\alpha\wedge\omega_{1})$ is satisfied. If the initial $SU(2)$-structure is hypo, we have $d\omega_{1}=0$, for all times $t$. So
$$(\frac{1}{2}\omega_{3}^{2})^{\cdot}=(\frac{1}{2}\omega_{1}^{2})^{\cdot}=\omega_{1}\wedge\dot{\omega}_{1}=\omega_{1}\wedge d\alpha=d(\alpha\wedge\omega_{1})$$
and we obtain the following $SU(2)$ and $SU(3)$-analogue of Theorem \ref{thmltex}.\\

\begin{corollary}\label{corltex}\textup{
There are no nontrivial longtime solutions for the hypo $SU(2)$ and $SU(3)$ embedding problem on compact manifolds.\\}
\end{corollary}
\rightline{$\Box$}


In view of Proposition \ref{propeea}, the following theorem yields solutions of the $G_{2}$ embedding problem.\\

\begin{theorem}\label{thmmain}\textup{
Let $\psi$ be a real analytic hypo $G_{2}$-structure on the compact manifold $M^{7}$. Then the intrinsic torsion flow
\begin{equation*}
\begin{cases}
\dot{A}_{t}=\mathcal{T}_{t}\circ A_{t}\\
A_{0}=\text{id}
\end{cases}
\end{equation*}
has a unique real analytic solution $A:(-\varepsilon,\varepsilon)\times M\rightarrow\text{End}(TM)$. Moreover, the solution $A_{t}$ is of the form
$$A_{t}=\sum_{k=0}^{\infty}\frac{t^{k}}{k!}A^{(k)}_{0},$$
where the series converges in the $C^{\infty}$-topology on $C^{\infty}(\text{End}(TM))$.}
\end{theorem}

\begin{beweis}
To apply Theorem \ref{thmckglobal} we have to show that the map
$$X:C^{\infty}(\text{Aut}(TM))\rightarrow C^{\infty}(\text{End}(TM))\;\;\;\text{ with }\;\;\;X\circ A:=\mathcal{T}(A\varphi)\circ A$$
is a real analytic first order differential operator in the sense of Definition \ref{defradiffop}. For this choose local coordinates $u:U\subset M\rightarrow\mathbb{R}^{7}$, for which $\varphi$ is real-analytic. These coordinates induce a local trivialization $(\pi,v)$ of the bundle $\pi:\text{End}(TM)\rightarrow M$ via
$$v(A):=\{a_{kl}\}_{k,l=1..7},\;\;\;\text{ where }\;\;\;A=\sum_{k,l=1}^{7}a_{kl}du_{k}\otimes\frac{\partial}{\partial u_{l}}\in\text{End}(TM).$$
For a fixed local section $A=\sum a_{kl}du_{k}\otimes\frac{\partial}{\partial u_{l}}:U\rightarrow\text{Aut}(TM)$ write
$$X\circ A=\mathcal{T}(A\varphi)\circ A=\sum_{a,b=1}^{7}f_{ab}du_{a}\otimes\frac{\partial}{\partial u_{b}}.$$
Now it suffices to find an expression
\setcounter{equation}{0}
\begin{equation}
f_{ab}=G_{ab}(u,a_{kl},\frac{\partial a_{kl}}{\partial u_{j}})
\end{equation}
for the coefficients $f_{ab}:U\rightarrow\mathbb{R}$, where $G_{ab}:D\subset \mathbb{R}^{7}\times \mathbb{R}^{49}\times\mathbb{R}^{343}\rightarrow\mathbb{R}$ is real analytic. The formula
\begin{equation*}
\nabla^{Ag}A\varphi=-\mathcal{T}(A\varphi)\lrcorner(A\psi)
\end{equation*}
shows that the intrinsic torsion is a first order invariant of the $G_{2}$-structure and hence we can find an expression of the form (1) that is actually polynomial in $a_{kl}$ and $\frac{\partial a_{kl}}{\partial u_{j}}$, and real analytic in $u$, since the initial structure is real analytic.\\
\end{beweis}


\begin{lemma}\label{lemtordiff}\textup{
Suppose $\psi$ is a $G_{2}$-structure on $M$ and $F\in\text{Diff}(M)$. Then the intrinsic torsion satisfies
$$\mathcal{T}(F^{\ast}\psi)=F^{\ast}\mathcal{T}(\psi)=F^{-1}_{\ast}\mathcal{T}(\psi)F_{\ast}.$$}
\end{lemma}

\begin{beweis}
By Koszul's formula we have $F_{\ast}(\nabla^{F^{\ast}g}_{X}Y)=\nabla^{g}_{F_{\ast}X}F_{\ast}Y$ and hence
$$(\nabla^{F^{\ast}g}_{X}F^{\ast}\varphi)=F^{\ast}(\nabla^{g}_{F_{\ast}X}\varphi).$$
Since $\nabla^{g}\varphi=-\mathcal{T}\lrcorner\psi$, we get
\begin{equation*}
\begin{split}
\mathcal{T}(F^{\ast}\psi)X\lrcorner F^{\ast}\psi&=-\nabla^{F^{\ast}g}_{X}F^{\ast}\varphi=-F^{\ast}(\nabla^{g}_{F_{\ast}X}\varphi)\\
&=F^{\ast}(T(\psi)F_{\ast}X\lrcorner\psi)=F_{\ast}^{-1}\mathcal{T}(\psi)F_{\ast}X\lrcorner F^{\ast}\psi\\
\end{split}
\end{equation*}
and the Lemma follows from the non-degeneracy of $F^{\ast}\psi$.\\
\end{beweis}


\begin{lemma}\label{lemintor0}\textup{
Suppose $\psi$ is a $G_{2}$-structure on $M^{7}=S^{1}\times..\times S^{1}\times M^{7-k}$, which is the hypo lift of some $SU(4-k)$-structure on $M^{7-k}$. Then the Ricci tensor $\text{Ric}$ of the metric $g=g(\psi)$ satisfies for each $S^{1}$-direction $\frac{\partial}{\partial\theta}$
$$L_{\frac{\partial}{\partial\theta}}\text{Ric}=\text{Ric}\frac{\partial}{\partial\theta}=d\theta\circ\text{Ric}=0.$$
The intrinsic torsion $\mathcal{T}$ satisfies
$$L_{\frac{\partial}{\partial\theta}}\mathcal{T}=\mathcal{T}\frac{\partial}{\partial\theta}=0$$
and $d\theta\circ\mathcal{T}=0$ if the structure is hypo.\\}
\end{lemma}

\begin{beweis}
If $\psi$ is the hypo lift of some structure on $M^{7-k}$, then $g=d\theta_{1}^{2}+..+d\theta_{k}^{2}+g_{7-k}$, for some metric $g_{7-k}$ on $M^{7-k}$. Hence the Ricci tensor satisfies $\text{Ric}\frac{\partial}{\partial\theta}=0$,
$$d\theta\circ\text{Ric}=g(\frac{\partial}{\partial\theta},\text{Ric})=g(\text{Ric}\frac{\partial}{\partial\theta},.)=0$$
and
$$L_{\frac{\partial}{\partial\theta}}\text{Ric}=\left.\frac{\partial}{\partial s}\right|_{s=0}\Phi^{\ast}_{s}\text{Ric}(g)
=\left.\frac{\partial}{\partial s}\right|_{s=0}\text{Ric}(\Phi^{\ast}_{s}g)
=\left.\frac{\partial}{\partial s}\right|_{s=0}\text{Ric}(g)=0.$$
Since $0=\nabla^{g}_{\frac{\partial}{\partial\theta}}\varphi=-\mathcal{T}\frac{\partial}{\partial\theta}\lrcorner\psi$, we get $\mathcal{T}\frac{\partial}{\partial\theta}=0$ and similarly $0=(L_{\frac{\partial}{\partial\theta}}\mathcal{T})\lrcorner\psi$ implies $L_{\frac{\partial}{\partial\theta}}\mathcal{T}=0$. If the structure is hypo, i.e. $\mathcal{T}$ is symmetric, we get in addition
$$d\theta\circ\mathcal{T}=g(\frac{\partial}{\partial\theta},\mathcal{T})=g(\mathcal{T}\frac{\partial}{\partial\theta},.)=0.$$
\end{beweis}


\begin{lemma}\label{lemapsi}\textup{
Suppose $\psi$ is a $G_{2}$-structure on $M^{7}=S^{1}\times..\times S^{1}\times M^{7-k}$, which is the hypo lift of some $SU(4-k)$-structure on $M^{7-k}$. If $A\in C^{\infty}(\text{Aut}(TM))$ satisfies
$$A\frac{\partial}{\partial\theta_{i}}=\frac{\partial}{\partial\theta_{i}},\;\;\;d\theta_{i}\circ A=d\theta_{i}\;\;\;\text{ and }\;\;\;
L_{\frac{\partial}{\partial\theta_{i}}}A=0,$$
then $A\psi$ is still the hypo lift of some $SU(4-k)$-structure.}
\end{lemma}

\begin{beweis}
By Lemma \ref{lemishl} we have $L_{\frac{\partial}{\partial\theta_{i}}}(A\psi)=0$ and
$$(Ag)(\frac{\partial}{\partial\theta_{i}},X)=g(\frac{\partial}{\partial\theta_{i}},A^{-1}X)=d\theta_{i} (A^{-1}X)=d\theta_{i}(X)
=g(\frac{\partial}{\partial\theta_{i}},X).$$
Now the Lemma follows from Lemma \ref{lemishl}.\\
\end{beweis}


We can now state the main result of this section,\\

\begin{theorem}\label{thmmain2}\textup{
Suppose $\psi$ is a real analytic hypo $G_{2}$-structure on $M=S^{1}\times..\times S^{1}\times M^{7-k}$, which is the hypo lift of some $SU(4-k)$-structure on $M^{7-k}$. Then the solution $A_{t}$ of the intrinsic torsion flow from Theorem \ref{thmmain} satisfies
$$A_{t}\frac{\partial}{\partial\theta_{i}}=\frac{\partial}{\partial\theta_{i}},\;\;\;d\theta_{i}\circ A_{t}=d\theta_{i}\;\;\;\text{ and }\;\;\;
L_{\frac{\partial}{\partial\theta_{i}}}A_{t}=0.$$
In particular, $A_{t}\psi$ is the hypo lift of some family of $SU(4-k)$-structures on $M^{7-k}$.}
\end{theorem}

\begin{beweis}
We apply Corollary \ref{corintcurv} with the following dictionary,
\begin{equation*}
\begin{array}{cl}
(1) &\mathcal{F}:=C^{\infty}(\text{End}(TM))\times C^{\infty}(\text{End}(TM))\\\\
(2) &\mathcal{U}:=C^{\infty}(\text{Aut}(TM))\times C^{\infty}(\text{End}(TM))\\\\
(3) &\mathcal{E}:=\{(B,\mathcal{T})\in\mathcal{F}\mid
0=L_{\frac{\partial}{\partial\theta_{i}}}B=L_{\frac{\partial}{\partial\theta_{i}}}\mathcal{T}\text{ and}\\ &\hspace{2.9cm}\,0=B\frac{\partial}{\partial\theta_{i}}=\mathcal{T}\frac{\partial}{\partial\theta_{i}}
=d\theta_{i}(B)=d\theta_{i}(\mathcal{T})\}\\\\
(4) &X:\mathcal{U}\rightarrow\mathcal{F}\text{ is defined w.r.t. the initial metric }g,\\
&\left.X\right|_{(A,\mathcal{T})}:=\big(\mathcal{T}\circ A,\text{Ric}(Ag)-\text{tr}(\mathcal{T})\mathcal{T}\big).\\\\
(5) &c(t):=(A_{t},\mathcal{T}_{t}).\\
\end{array}
\end{equation*}
Note that $\mathcal{U}\subset\mathcal{F}$ is open by Example \ref{example}, and that $X$ is smooth and $\mathcal{E}\subset\mathcal{F}$ is closed, since differential operators are smooth by Example 3.6.6. in \cite{ham}. By Proposition \ref{propeea}, Theorem \ref{thmgt} and the definition of $A_{t}$, the curve $c(t)$ is an integral curve of the vector field $X$.
From Lemma \ref{lemintor0} we get $c(0)=(\text{id},\mathcal{T}_{0})\in\mathcal{E}_{f}$, where $f:=(\text{id},0)\in\mathcal{F}$. Now it suffices to show that $X$ is tangent to $\mathcal{U}\cap\mathcal{E}_{f}$, i.e.
$$X_{\mid \mathcal{U}\cap\mathcal{E}_{f}}:\mathcal{U}\cap\mathcal{E}_{f}\rightarrow\mathcal{E}.$$
For $(A=\text{id}+B,\mathcal{T})\in\mathcal{U}\cap\mathcal{E}_{f}$ we have
$$A\frac{\partial}{\partial\theta_{i}}=\frac{\partial}{\partial\theta_{i}},\;\;\;d\theta_{i}\circ A=d\theta_{i}\;\;\;\text{ and }\;\;\;
L_{\frac{\partial}{\partial\theta_{i}}}A=0.$$
By Lemma \ref{lemapsi} we see that $A\psi$ is still the hypo lift of some $SU(4-k)$-structure and Lemma \ref{lemintor0} yields
$$L_{\frac{\partial}{\partial\theta_{i}}}\text{Ric}(Ag)=\text{Ric}(Ag)\frac{\partial}{\partial\theta_{i}}=d\theta_{i}\circ\text{Ric}(Ag)=0.$$
Now we can easily verify that $X(A,\mathcal{T})\in\mathcal{E}$,
\begin{enumerate}
\item[$\bullet$] $L_{\frac{\partial}{\partial\theta_{i}}}(\mathcal{T}\circ A)=0$ and $L_{\frac{\partial}{\partial\theta_{i}}}(\text{Ric}(Ag)-\text{tr}(\mathcal{T})\mathcal{T})=0,$\\
\item[$\bullet$] $\mathcal{T}\circ A\frac{\partial}{\partial\theta_{i}}=0$
and $(\text{Ric}(Ag)-\text{tr}(\mathcal{T})\mathcal{T})\frac{\partial}{\partial\theta_{i}}=0$,\\
\item[$\bullet$] $d\theta_{i}(\mathcal{T}\circ A)=0$ and
$d\theta_{i}(\text{Ric}(Ag)-\text{tr}(\mathcal{T})\mathcal{T})=0$\\
\end{enumerate}
and the Theorem follows.\\
\end{beweis}


\begin{remark}\label{remmain2}\textup{
The property $L_{\frac{\partial}{\partial\theta}}A_{t}=0$ from Theorem \ref{thmmain} is a consequence of the diffeomorphism invariance of the evolution equation $\dot{A}_{t}=\mathcal{T}_{t}\circ A_{t}$. In fact, Lemma \ref{lemtordiff} shows that $B_{t}:=\Phi_{s}^{\ast}A_{t}$ also solves $\dot{A}_{t}=\mathcal{T}_{t}\circ A_{t}$, where $\Phi_{s}$ is the flow of $\frac{\partial}{\partial\theta}$. Since $\Phi_{s}$ is real analytic, the uniqueness part of Theorem \ref{thmmain} yields $A_{t}=\Phi_{s}^{\ast}A_{t}$, i.e. $L_{\frac{\partial}{\partial\theta}}A_{t}=0$.\\}
\end{remark}


We can now solve the embedding problem for real analytic hypo $SU(4-k)$-structures on $M^{7-k}$ by reducing it to the embedding problem for real analytic hypo $G_{2}$-structures on $M=S^{1}\times..\times S^{1}\times M^{7-k}$. Namely, the hypo lift of the initial $SU(4-k)$-structure yields a real analytic hypo $G_{2}$-structures on $M$. Theorem \ref{thmmain} yields a solution $A_{t}$ of the intrinsic torsion flow. By Theorem \ref{thmmain2} the family of $G_{2}$-structures $\psi_{t}=A_{t}\psi$ is still the hypo lift of some family of $SU(4-k)$-structures. Now Lemma \ref{lemhlee} proves that the family of $SU(4-k)$-structures is a solution of the embedding problem.\\

\begin{corollary}\label{cormain2}\textup{
For any real analytic hypo $SU(2)$, $SU(3)$ and $G_{2}$-structure on a compact manifold, the embedding problem admits a unique real analytic solution. Moreover, the solution can be described by a family of gauge deformations
$$A_{t}=\sum_{k=0}^{\infty}\frac{t^{k}}{k!}A^{(k)}_{0},$$
where the series converges in the $C^{\infty}$-topology on $C^{\infty}(\text{End}(TM))$.\\}
\end{corollary}


\setcounter{chapter}{4}
\setcounter{theorem}{0}
\vspace{1cm}
\begin{center}\label{su2}{\textsc{\textbf{Appendix: $SU(2)$-Structures in Dimension Five}}}\end{center}
\vspace{1cm}

Usually a $SU(2)$-structure on a five dimensional manifold is described by a quadruplet of forms $(\alpha,\omega_{1},\omega_{2},\omega_{3})$, cf. for instance \cite{cosa}. There is an alternative to the usual definition, which is justified by the last equation in the next Lemma.\\

\begin{lemma}\label{lemdefsu2}\textup{
\begin{equation*}
\begin{split}
\text{Iso}_{GL(5)}(\alpha_{0})&=\{\begin{pmatrix}1 &0\\ x &A\end{pmatrix}\mid A\in GL(4)\text{ and } x\in\mathbb{R}^{4}\}.\\
\text{Iso}_{GL(5)}(\omega_{1})&=\{\begin{pmatrix}\lambda &y^{T}\\ 0 &A\end{pmatrix}\mid A\in \text{Sp}(4,\mathbb{R}), y\in\mathbb{R}^{4}\text{ and }\lambda\neq 0\}.\\
\text{Iso}_{GL(5)}(\alpha_{0},\omega_{1},\omega_{2},\omega_{3})&=\begin{pmatrix}1&0\\0&SU(2)\end{pmatrix}.\\
\text{Iso}_{GL^{+}(5)}(\omega_{1},\rho_{2},\rho_{3})&=\begin{pmatrix}1&0\\0&SU(2)\end{pmatrix}.\\
\end{split}
\end{equation*}}
\end{lemma}

\begin{beweis}
Write $B\in GL(5)$ as
$$B=\begin{pmatrix}\lambda &y^{T}\\ x &A\end{pmatrix},$$
where $\lambda\in\mathbb{R}$, $x,y\in\mathbb{R}^{4}$ and $A\in\mathfrak{gl}(4)$. Then $\alpha(Be_{1})=\lambda$ and $\alpha(Be_{j})=y^{T}e_{j}$, for $j\in\{2,..,5\}$. Hence the stabilizer of the $1$-form $\alpha_{0}:=e^{1}\in\Lambda^{1}\mathbb{R}^{5\ast}$ has the above form.\\
For $B\in\text{Iso}_{GL(5)}(\omega_{1})$ and $i,j\in\{2,..,5\}$ we get $\omega_{1}(e_{i},e_{j})=\omega_{1}(Be_{i},Be_{j})=\omega_{1}(Ae_{i},Ae_{j})$, i.e. $A\in\text{Sp}(4,\mathbb{R})$. This yields
$$0=\omega_{1}(Be_{1},Be_{j})=\omega_{1}(\lambda e_{1}+x,(y^{T}e_{j})e_{1}+Ae_{j})=\omega_{1}(x,Ae_{j})=\omega_{1}(A^{-1}x,e_{j})$$
and the non-degeneracy of $\omega_{1}$, as a form on $\mathbb{R}^{4}$, implies $x=0$ and proves the second equation of the lemma.\\
Now the third equation follows, since $\omega_{2}=\text{Re}(\Phi_{0})$ and $\omega_{3}=\text{Im}(\Phi_{0})$, where $\Phi_{0}=(e^{2}+ie^{3})\wedge(e^{4}+ie^{5})$, and $SU(2)=\text{Sp}(4,\mathbb{R})\cap SL(2,\mathbb{C})$.\\
To obtain the last equation, we compute for $B=\begin{pmatrix}\lambda &y^{T}\\ 0 &A\end{pmatrix}\in\text{Iso}_{GL(5)}(\omega_{1})\cap\text{Iso}_{GL^{+}(5)}(\alpha_{0}\wedge\omega_{2})$ and $i,j\in\{2,..,5\}$
\begin{equation*}
\begin{split}
\omega_{2}(e_{i},e_{j})&=(\alpha_{0}\wedge\omega_{2})(e_{1},e_{i},e_{j})=(\alpha_{0}\wedge\omega_{2})(Be_{1},Be_{i},Be_{j})\\
&=(\alpha_{0}\wedge\omega_{2})(\lambda e_{1},(y^{T}e_{i})e_{1}+Ae_{i},(y^{T}e_{j})e_{1}+Ae_{j})\\
&=(\alpha_{0}\wedge\omega_{2})(\lambda e_{1},Ae_{i},Ae_{j})\\
&=\lambda\omega_{2}(Ae_{i},Ae_{j}).\\
\end{split}
\end{equation*}
Since the volume element $\varepsilon_{0}=e^{2345}$ on $\mathbb{R}^{4}$ satisfies
$$\varepsilon_{0}=\frac{1}{2}\omega_{1}^{2}=\frac{1}{2}\omega_{2}^{2}=\frac{1}{2}\omega_{3}^{2},$$
we obtain from $A\in\text{Sp}(4,\mathbb{R})=\text{Iso}_{GL(4)}(\omega_{1})$
$$\text{det}(A)\varepsilon_{0}=A^{-1}\varepsilon_{0}=A^{-1}\frac{1}{2}\omega_{1}^{2}=\varepsilon_{0},$$
i.e. $\text{det}(A)=1$. Now $A^{-1}\omega_{2}=\lambda^{-1}\omega_{2}$ yields
$$\varepsilon_{0}=A^{-1}\frac{1}{2}\omega_{1}^{2}=\lambda^{-2}\varepsilon_{0}$$
and since $B\in GL^{+}(5)$, we get $\lambda=1$. Similarly we get $A\omega_{3}=\omega_{3}$, which yields $A\in SU(2)$. Now
\begin{equation*}
\begin{split}
\alpha_{0}\wedge\omega_{2}&=B^{-1}(\alpha_{0}\wedge\omega_{2})=B^{-1}\alpha_{0}\wedge B^{-1}\omega_{2}\\
&=B^{-1}\alpha_{0}\wedge A^{-1}\omega_{2},\;\;\;\text{ since }e_{1}\lrcorner\omega_{2}=0\\
&=(\alpha_{0}(Be_{1})e^{1}+\sum_{j=2}^{5}\alpha_{0}(Be_{j})e^{j})\wedge\omega_{2}\\
&=(\alpha_{0}+\sum_{j=2}^{5}y_{j}e^{j})\wedge\omega_{2}
\end{split}
\end{equation*}
yields $\sum_{j=2}^{5}y_{j}e^{j}\wedge\omega_{2}=0$, i.e. $y=0$.\\
\end{beweis}


Since the $GL^{+}(5)$ stabilizer of the triple $(\omega_{1},\rho_{2},\rho_{3})$ is equal to $\{1\}\times SU(2)$, we expect that, after fixing an orientation for $\mathbb{R}^{5}$, we can reconstruct the forms $\alpha_{0}$, $\omega_{2}$ and $\omega_{3}$ solely from the triple $(\omega_{1},\rho_{2},\rho_{3})$. The first step is to reconstruct the volume element $\varepsilon_{0}$. Then the forms $\alpha_{0}$, $\omega_{2}$ and $\omega_{3}$, as well as the metric $g_{0}$, can be obtained from the formulas in Example \ref{exsu2}.\\


\begin{lemma}\label{lemepssu2}\textup{
After choosing an orientation for $V:=\mathbb{R}^{5}$, there is a homomorphism
$$\varepsilon:\Lambda^{2}V^{\ast}\oplus\Lambda^{3}V^{\ast}\oplus\Lambda^{3}V^{\ast}\rightarrow\Lambda^{5}V^{\ast}\oplus i\Lambda^{5}V^{\ast}$$
of $GL^{+}(5)$-modules, such that for the model tensors and the canonical orientation $[\varepsilon_{0}]$ of $\mathbb{R}^{5}$
$$\varepsilon(\omega_{1},\rho_{2},\rho_{3})=\varepsilon_{0}\in\Lambda^{5}V^{\ast}\subset\Lambda^{5}V^{\ast}\oplus i\Lambda^{5}V^{\ast}.$$}
\end{lemma}

\begin{beweis}
Given an orientation $[\varepsilon_{+}]$ for $V$, represented by an element $\varepsilon_{+}\in\Lambda^{5}V^{\ast}$, we can define a $GL^{+}(5)$-equivariant map
$$\sqrt[4]{\;\;}:\Lambda^{5}V^{\ast}\otimes\Lambda^{5}V^{\ast}\otimes\Lambda^{5}V^{\ast}\otimes\Lambda^{5}V^{\ast}\rightarrow\Lambda^{5}V^{\ast}\oplus i\Lambda^{5}V^{\ast}.$$
Now consider the $GL(5)$-equivariant map
$$K:\Lambda^{2}V^{\ast}\oplus\Lambda^{3}V^{\ast}\oplus\Lambda^{3}V^{\ast}\rightarrow (V^{\ast}\otimes V)\otimes(V^{\ast}\otimes V)
\otimes\Lambda^{5}V^{\ast}\otimes\Lambda^{5}V^{\ast}$$
defined by
$$K(\omega_{1},\rho_{2},\rho_{3})(x,a,y,b):=\big(\rho_{2}\wedge a\wedge b\big)
\otimes\big(\rho_{3}\wedge(x\lrcorner\omega_{1})\wedge(y\lrcorner\omega_{1})\big),$$
where $x,y\in V$ and $a,b\in V^{\ast}$. For the model tensors $\omega_{1},\rho_{2},\rho_{3}$ let $K_{0}:=K(\omega_{1},\rho_{2},\rho_{3})$. Then we compute
\begin{equation*}
\begin{split}
K_{0}(x,a,y,b)=(a_{5}b_{3}-a_{3}b_{5}+a_{2}b_{4}-a_{4}b_{2})(-x_{3}y_{4}+x_{4}y_{3}-x_{2}y_{5}+x_{5}y_{2})\otimes\varepsilon_{0}^{2}.
\end{split}
\end{equation*}
Taking the trace of the first factor $V^{\ast}\otimes V$, we obtain a map
$$L=\text{tr}(K):\Lambda^{2}V^{\ast}\oplus\Lambda^{3}V^{\ast}\oplus\Lambda^{3}V^{\ast}\rightarrow (V^{\ast}\otimes V)
\otimes\Lambda^{5}V^{\ast}\otimes\Lambda^{5}V^{\ast}$$
and for the model tensors we obtain
\begin{equation*}
\begin{split}
L_{0}(y,b):=\text{tr}(K_{0})(y,b)&=(-b_{4}y_{5}+b_{5}y_{4}-b_{2}y_{3}+b_{3}y_{2})\otimes\varepsilon_{0}^{2}.
\end{split}
\end{equation*}
Identifying $V^{\ast}\otimes V=\text{Hom}(V,V)$, we define
$$L^{2}:\Lambda^{2}V^{\ast}\oplus\Lambda^{3}V^{\ast}\oplus\Lambda^{3}V^{\ast}\rightarrow (V^{\ast}\otimes V)
\otimes(\Lambda^{5}V^{\ast})^{4}$$
and so
$$L_{0}^{2}=\begin{pmatrix}0&0\\0&-\text{id}_{\mathbb{R}^{4}}\end{pmatrix}\otimes\varepsilon_{0}^{4}.$$
Taking again the trace, we obtain a map
$$\text{tr}(L^{2}):\Lambda^{2}V^{\ast}\oplus\Lambda^{3}V^{\ast}\oplus\Lambda^{3}V^{\ast}\rightarrow(\Lambda^{5}V^{\ast})^{4}$$
with $\text{tr}(L_{0}^{2})=-4\varepsilon_{0}^{4}$. Hence
$$\varepsilon:=\sqrt[4]{-\frac{1}{4}\text{tr}(L^{2})}
:\Lambda^{2}V^{\ast}\oplus\Lambda^{3}V^{\ast}\oplus\Lambda^{3}V^{\ast}\rightarrow\Lambda^{5}V^{\ast}\oplus i\Lambda^{5}V^{\ast}$$
is the desired equivariant map.\\
\end{beweis}


\vspace{1cm}

\renewcommand\bibname{References}



\begin{thebibliography}{00}

\bibitem{baer}C.Bär, P.Gauduchon, A.Moroianu:
\textit{Generalized cylinders in semi Riemannian and spin
geometry}, Mathematische Zeitschrift 249 545 (2005).

\bibitem{bes}A.L. Besse:
\textit{Einstein Manifolds}, Springer Verlag, Berlin-Heidelberg-New York, 1987.

\bibitem{bry}R.Bryant:
\textit{Nonembedding and Nonextension Results in Special Holonomy}, In Proceedings
of the August 2006 Madrid conference in honor of Nigel Hitchin's 60th Birthday. Oxford University Press.

\bibitem{bryeds}R.Bryant et al.:
\textit{Exterior Differntial Systems}, Mathematical Sciences Research Institute Publications, Springer -Verlag, 1991.


\bibitem{car}H.Cartan:
\textit{Elementare Theorie der analytischen Funktionen einer oder mehrerer komplexen Ver\"{a}nderlichen}, 1966, Bibliographisches Institut AG, Mannheim.

\bibitem{co}D.Conti:
\textit{Embedding into Manifolds with Torsion}, arXiv:0812.4186, 2009.

\bibitem{cosa}D.Conti, S.Salamon:
\textit{Generalized Killing Spinors in Dimension 5}, Trans. Amer. Math. Soc., 359(11):5319-5343, 2007.

\bibitem{cor}V.Cort\'{e}s, T.Leistner, L.Sch\"{a}fer, F.Schulte-Hengesbach:
\textit{Half-flat Structures and Special Holonomy}, arXiv:0907.1222, 2009.


\bibitem{fer1}M. Fernandez, S.Ivanov, V.Munoz,
L.Ugarte:
\textit{Nearly hypo structures and compact nearly
Kähler 6-manifolds with conical singularities}, to appear in the Journal of the London. Math. Soc., arXiv:math/0602160.

\bibitem{fer2}M. Fernandez, A. Tomassini, L. Ugarte, R. Villacampa:
\textit{Balanced Hermitian Metrics from $SU(2)$-Structures}, Journal of Mathematical Physics, Volume 50, Issue 3, 2009.

\bibitem{friv}T. Friedrich, S. Ivanov:
\textit{Parallel Spinors and Connections with
Skew-Symmetric Torsion in String Theory}, AsianJ.Math.6:303-336, 2002.


\bibitem{ham}R.S. Hamilton:
\textit{The Inverse Function Theorem of Nash and Moser}, Bull. Amer. Math. Soc. 7,
1982, pages 65-222.

\bibitem{hit1}N. Hitchin:
\textit{Stable forms and special metrics}, Global Differential Geometry: \textit{The
Mathematical Legacy of Alfred Gray}, volume 288 of Contemp. Math.,
pages 70-89. American Math. Soc., 2001.


\bibitem{sal}S. Salamon:
\textit{Riemannian Geometry and Holonomy Groups}, Pitman Research Notes in Mathematics 201, Longman, 1989.


\end{thebibliography}
\end{document}